\documentclass[a4paper,11pt]{article}
\newif\ifpdf
\ifpdf
\fi
\usepackage{amssymb}
\usepackage{enumerate}

\newtheorem{tw}{Theorem}[subsection]

\newtheorem{lm}[tw]{Lemma}

\newtheorem{wn}[tw]{Corollary}
\newtheorem{stw}[tw]{Proposition}
\newenvironment{dow}{\it Proof.\rm}{\hfill $\Box$}
\newcommand{\Rightsx}{\mathop{\Rightarrow}_{s,x}}

\newcommand{\nsubsection}{\setcounter{equation}{0}\subsection}

\begin{document}

\title {Reflected BSDEs and continuous solutions of parabolic obstacle
problem for semilinear PDEs in divergence form}
\author {Tomasz Klimsiak}
 \maketitle
\begin{abstract}
We consider the Cauchy problem for semilinear parabolic equation
in divergence form with obstacle.  We show that under natural
conditions on the right-hand side of the equation and mild
conditions on the obstacle a unique continuous solution of the
problem admits a stochastic representation in terms of reflected
backward stochastic differential equations. We derive also some
regularity properties of solutions  and prove useful approximation
results.
\end{abstract}

\footnotetext{T. Klimsiak: Faculty of Mathematics and Computer
Science, Nicolaus Copernicus University, tomas@mat.uni.torun.pl}


\footnotetext{{\em Key words or phrases:} Backward stochastic
differential equation, Semilinear parabolic partial differential
equation, Divergence form operator, Obstacle problem, Weak
solution.}

\nsubsection{Introduction}

In the present paper we are interested in stochastic representation
of solutions of the Cauchy problem
for semilinear parabolic equation in divergence form with obstacle.
Let $a:Q_T\equiv[0,T]\times\mathbb{R}^{d}
\rightarrow\mathbb{R}^{d}\otimes\mathbb{R}^{d}$ be a measurable,
symmetric matrix valued  function such that
\begin{equation}\label{eq1.1}
\lambda|\xi|^2\leq\sum^d_{i,j=1}a_{ij}(t,x)\xi_i\xi_j\le\Lambda|\xi|^2,
\quad \xi\in\mathbb{R}^d
\end{equation}
for some $0<\lambda\le\Lambda$ and let $A_t$ be a linear operator of
the form
\begin{equation}\label{eq1.2}
A_{t}=\frac{1}{2}\sum_{i,j=1}^d\frac{\partial}{\partial
x_{i}}(a^{ij}\frac{\partial}{\partial x_{j}}).
\end{equation}
Roughly speaking the problem consist in finding
$u:Q_T
\rightarrow\mathbb{R}$ such that for given
$\varphi:\mathbb{R}^d\rightarrow\mathbb{R}$,
$h:Q_T\rightarrow\mathbb{R}$, $f:Q_T\times\mathbb{R}
\times\mathbb{R}^d\rightarrow\mathbb{R}$,
\begin{equation}
\label{eq1.3} \left\{
\begin{array}{ll} \min(u-h,-\frac{\partial u}{\partial t}
-A_{t}u-f_u)=0 & \mbox{in }Q_T,\\
u(T)=\varphi & \mbox{on }\mathbb{R}^d,
\end{array}
\right.
\end{equation}
where $f_u=f(\cdot,\cdot,u,\sigma\nabla u)$ and
$\sigma\sigma^*=a$, i.e. $u$ satisfies the prescribed terminal
condition, takes values above a given obstacle $h$,  satisfies
inequality $\frac{\partial u}{\partial t}+A_tu\leq -f_u$ in $Q_T$
and equation $\frac{\partial u}{\partial t}+A_tu=-f_u$ on the set
$\{u>h\}$.

The obstacle problem (\ref{eq1.3}) has been studied intensively by
many authors.  Subject to regularity of the data $\varphi,f,h$ and
coefficients of $A_{t}$, viscosity solutions (see \cite{EKPPQ}) or
solutions of variational inequalities associated with
(\ref{eq1.3}) are considered. In the latter case  one can consider
weak solutions (see \cite{BensoussanLions,Lions,Mignot}) or strong
solutions (see \cite{BensoussanLions,Brezis,Charrier,Donati}).

In the present paper by a solution of (\ref{eq1.3}) we understand
a pair $(u,\mu)$ consisting of a measurable function
$u:Q_T\rightarrow\mathbb{R}$ having some regularity properties and
a Radon measure $\mu$ on $Q_T$ such that
\begin{equation}
\label{eq1.4} \frac{\partial u}{\partial t }+A_tu=-f_u-\mu,\quad
u(T)=\varphi,\quad u\geq h,\quad \int_{Q_T}(u-h)\,d\mu=0
\end{equation}
(see Section \ref{ssec2.2} for details). We adopt the above
definition for three reasons. Firstly, it may be viewed as an
analogue of the definition of the obstacle problem for elliptic
equations (see \cite{Kinderlehrer,Leone}). It is worth pointing out,
however, that contrary to the case of elliptic equations, it is not
obvious how solution of a parabolic variational inequality
associated with (\ref{eq1.3}) is related to the solution  in the
sense of (\ref{eq1.4}). Secondly, since in many cases we are able to
prove some additional information on $\mu$, using (\ref{eq1.4})
instead of variational formulation gives more information on
solutions of (\ref{eq1.3}). Finally, definition (\ref{eq1.4}) is
well suited with our main purpose which consists in providing
stochastic representation of solutions of the obstacle problem.

In the case where $A_t$ is a non-divergent operator of the form
\[
A_{t}=\frac{1}{2}\sum_{i,j=1}^da^{ij} \frac{\partial^{2}}{\partial
x_{i}\partial x_{j}}+\sum_{i=1}^db_{i}\frac{\partial}{\partial
x_{i}}\,,
\]
problem (\ref{eq1.3}) has been investigated carefully in
\cite{EKPPQ} by using probabilistic methods. Let $X^{s,x}$ be a
solution of the It\^o equation
\[
dX^{s,x}_{t}=\sigma(t,X^{s,x}_t)dW_t+b(t,X^{s,x}_t)\,dt, \quad
X^{s,x}_s=x\quad (\sigma\sigma^*=a)
\]
associated with $A_t$. In \cite{EKPPQ} it is proved, that under
suitable assumptions on $a,b$ and the data $\varphi,f,h$, for each
$(s,x)\in Q_T$  there exists a unique solution
$(Y^{s,x},Z^{s,x},K^{s,x})$ of reflected backward stochastic
differential equation with forward driving process $X^{s,x}$,
terminal condition $\varphi(X^{s,x}_T)$, coefficient $f$ and
obstacle $h(\cdot,X^{s,x}_{\cdot})$ (RBSDE$(\varphi,f,h)$ for
short), and moreover, $u$ defined by the formula
$u(s,x)=Y^{s,x}_s$, $(s,x)\in Q_T$ is a unique viscosity solution
of (\ref{eq1.3}) in the class of functions satisfying the
polynomial growth condition. In the present paper we give a
representation similar to that proved in \cite{EKPPQ} for weak
solutions of (\ref{eq1.4}) with $A_t$ defined by (\ref{eq1.2}).

In the paper we assume that
\begin{enumerate}
\item[(H1)]$\varphi\in\mathbb{L}_2^{loc}(\mathbb{R}^{d})$,
$h\in\mathbb{L}_2^{loc}(Q_{T})$, 
\item[(H2)]$f:[0,T]\times\mathbb{R}^{d}\times\mathbb{R}\times\mathbb{R}^{d}
\rightarrow\mathbb{R}$ is a measurable function satisfying the
following conditions:
\begin{enumerate}[{a)}]
\item there is $L>0$ such that
$|f(t,x,y_{1},z_{1})-f(t,x,y_{2},z_{2})|\leq
L(|y_{1}-y_{2}|+|z_{1}-z_{2}|)$ for all
$(t,x)\in[0,T]\times\mathbb{R}^{d}$, $y_{1},y_{2}\in\mathbb{R}$
and $z_{1}, z_{2}\in\mathbb{R}^{d}$,
\item there exist $M>0$, $g\in\mathbb{L}^{loc}_{2}(Q_{T})$ such that
$|f(t,x,y,z)|\leq g(t,x)+M(|y|+|z|)$ for all
$(t,x,y,z)\in[0,T]\times\mathbb{R}^{d}\times\mathbb{R}\times\mathbb{R}^{d}$,
\end{enumerate}
\item[(H3)]$\varphi(x)\geq h(T,x)$ for all $x\in\mathbb{R}^{d}, \, h\in C(Q_{T})$
\end{enumerate}
(definitions of various function spaces used in the paper are
given at the end of the section).

We prove that under (\ref{eq1.1}) and (H1)--(H3) the obstacle
problem (\ref{eq1.4}) has at most one solution such that $u\in
C(\breve{Q}_{T})\cap W^{0,1}_{2,\varrho}(Q_T)$ for some $\varrho$
of the form $\varrho(x)=(1+|x|^2)^{-\alpha}$, $x\in\mathbb{R}^d$,
where $\alpha\ge0$. From our existence results it follows in
particular that if, in addition,
$\varphi\in\mathbb{L}_{2,\varrho}(\mathbb{R}^d)$,
$g\in\mathbb{L}_{2,\varrho}(Q_T)\cap\mathbb{L}_{p,q,\varrho}(Q_T)$,
$h\in C(Q_T)\cap\mathbb{L}_{2,\varrho}(Q_T)$ for some $\varrho$ as
above and $p,q\in(2,\infty]$ such that $(2/q)+(d/p)<1$, and $h$
satisfies the polynomial growth condition, then (\ref{eq1.4}) has
a solution $(u,\mu)$ such that $u\in C([0,T)\times
\mathbb{R}^{d})\cap W_{2,\varrho}^{0,1}(Q_{T})$. Secondly, for
each $(s,x)\in[0,T)\times\mathbb{R}^d$ we have
\begin{equation}
\label{eq1.5} (u(t,X_{t}),\sigma\nabla
u(t,X_{t}))=(Y^{s,x}_{t},Z^{s,x}_{t}),\quad t\in[s,T], \quad
P_{s,x}\mbox{-}a.s.,
\end{equation}
where $(X,P_{s,x})$ is a Markov process associated with $A_t$ (see
\cite{A.Roz.DIFF,Stroock.DIFF}) and $Y^{s,x},Z^{s,x}$ are the
first two components of a solution $(Y^{s,x},Z^{s,x},K^{s,x})$ of
RBSDE$(\varphi,f,h)$ with forward driving process $X$. In
particular, it follows that
\begin{equation}
\label{eq1.6} u(s,x)=Y^{s,x}_s,\quad (s,x)\in Q_T,
\end{equation}
which may be viewed as a generalization of the Feynman-Kac
formula. We show also that
\begin{equation}
\label{eq1.7} E_{s,x}\int_{s}^{T}\xi(t,X_t)\,dK^{s,x}_t
=\int_s^T\!\!\int_{\mathbb{R}^d}\xi(t,y)p(s,x,t,y)\,d\mu(t,y)
\end{equation}
for all $\xi\in C_{b}(Q_{T})$, where $p$ stands for the transition
density function of $(X,P_{s,x})$ (or, equivalently, $p$ is the
fundamental solution for $A_t$), which provides an additional
information on the process $K^{s,x}$ and solution $(u,\mu)$ of
(\ref{eq1.4}). For instance, it follows from (\ref{eq1.7}) that in
the linear case the solution of (\ref{eq1.3}) admits the
representation
\begin{eqnarray*}
&&u(s,x)=\int_{\mathbb{R}^d}\varphi(y)p(s,x,T,y)\,dy
+\int_{Q_T}f(t,y)p(s,x,t,y)\,dy\nonumber\\
&&\qquad\qquad+\int_{Q_T}p(s,x,t,y)\,d\mu(t,y),
\quad(s,x)\in[0,T)\times\mathbb{R}^d,
\end{eqnarray*}
which, up to our knowledge, is new (for parabolic problems).
Moreover, we show that $\mu$ is absolutely continuous with respect
to the Lebesgue measure $\lambda$ and $d\mu=r\,d\lambda$ if and
only if
\[
K^{s,x}_{t}=\int_{s}^{t}r(\theta,X^{s,x}_{\theta})\,d\theta,\quad
t\in[s,T].
\]

Let us remark also that the first component $u$ of a solution of
(\ref{eq1.4}) coincides with the solution of (\ref{eq1.3}) in the
variatonal sense.

Our conditions on  $\varphi,g$ and $h$ are similar to that used in
the theory of variational inequalities and seems to be close to
the best possible. As for $g$, in fact we prove existence and
uniqueness of solutions of (\ref{eq1.4}) and the representation
(\ref{eq1.5}) under the assumption that
$g\in\mathbb{L}_{2,\varrho}(Q_T)$ and
\begin{equation}
\label{eq1.8} E_{s,x}\int^T_s|g(t,X_t)|^2\,dt
\end{equation}
is bounded uniformly in $(s,x)\in K$ for every compact subset $K$
of $[0,T)\times\mathbb{R}^d$. We show also that if $\varphi\in
\mathbb{L}_{2,\varrho}(\mathbb{R}^{d})$, $g\in
\mathbb{L}_{2,\varrho}(Q_{T})$ and $h\in C(Q_T)$ satisfies the
polynomial growth condition, then there is a version of the
minimal weak solution of (\ref{eq1.3}) in the variational sense
such that if (\ref{eq1.8}) is finite for some fixed
$(s,x)\in[0,T)\times\mathbb{R}^d$, then  there exists a unique
solution $(Y^{s,x},Z^{s,x},K^{s,x})$ of RBSDE$(\varphi,f,h)$ and
(\ref{eq1.5}) holds true. Thus, since (\ref{eq1.8}) is finite for
a.e. $(s,x)\in[0,T)\times\mathbb{R}^d$ if
$g\in\mathbb{L}_{2,\varrho}(Q_{T})$, (\ref{eq1.5}) holds for a.e.
$(s,x)\in[0,T)\times\mathbb{R}^d$ if
$g\in\mathbb{L}_{2,\varrho}(Q_{T})$. What is more important, it
follows from our result that for each $(s,x)$ such that
(\ref{eq1.8}) is finite we get a probabilistic formula
(\ref{eq1.6})  for the minimal weak solution of the variational
inequality associated with (\ref{eq1.3}).

In case $\varrho=1$ existence of a solution of (\ref{eq1.4}) and
representation (\ref{eq1.5}) is proved by the method of stochastic
penalization used earlier in \cite{EKPPQ}. For $\varrho<1$ in
proofs of these results we use ideas from \cite{A.Roz.BSDE}. In
both cases from our proofs it follows that if $(u,\mu) $ is a
solution of (\ref{eq1.4}), $u_n$ is a solution of the Cauchy
problem
\[
(\frac{\partial}{\partial t}+A_t)
u_{n}=-f_{u_{n}}-n(u_{n}-h)^{-},\quad u_{n}(T)=\varphi
\]
and $\mu_n$ is a measure on $Q_T$ such that
$d\mu_n=n(u_n-h)^-\,d\lambda$  then $u_n\uparrow u$ uniformly in
compact subsets of $[0,T)\times\mathbb{R}^d$ and  in
$W^{0,1}_{2,\varrho}(Q_{T})\cap
C([0,T],\mathbb{L}_{2,\varrho}(Q_{T}))$ if $\varrho<1$, and
locally in the latter space if $\varrho=1$. In particular,
differently from the theory of variational inequalities, we obtain
strong convergence in $\mathbb{L}_{2,\varrho}(Q_{T})$ of gradients
of $u_n$'s to the gradient of $u$. Moreover, from the proofs it
follows that $\{\mu_n\}$ converges weakly to $\mu$  and strongly
in the space dual to $W^{1,1}_{2,\varrho}(Q_T)$, and for each
$(s,x)\in[0,T)\times\mathbb{R}^d$ the measures $\nu_n$ defined by
the relation $d\nu_n/d\mu_n=p(s,x,\cdot,\cdot)$ converge weakly to
the measure $\nu$ such that $d\nu/d\mu=p(s,x,\cdot,\cdot)$. These
results allow us to deduce some properties of $\mu$ from
properties of the sequence $\{\mu_n\}$.

In the paper we will use the following notation.

$Q_T=[0,T]\times\mathbb{R}^d$, $\breve
Q_T=(0,T)\times\mathbb{R}^d$. For $E\subset Q_T$ we write
$E_{t}=\{x\in\mathbb{R}^{d};(t,x)\in E\}$.
$B(0,r)=\{x\in\mathbb{R}^d:|x|<r\}$, $x^+=\max(x,0)$,
$x^-=\max(-x,0)$. $\nabla=(\frac{\partial}{\partial x_1},
\dots,\frac{\partial}{\partial x_d})$. By $\lambda$ we denote the
Lebesgue measure.

$\mathbb{L}_{p}(\mathbb{R}^d)$ is the Banach space of measurable
function $u$ on $\mathbb{R}^d$ having the finite norm
$\|u\|_p=(\int_{\mathbb{R}^d}|u(x)|^p\,dx)^{1/p}$.
$\mathbb{L}_{p,q}(Q_T)$ is the Banach space of measurable
functions on $Q_T$ having the finite norm
$\|u\|_{p,q,T}=(\int_0^T(\int_{\mathbb{R}^d}
|u(t,x)|^p\,dx)^{p/q}\,dt)^{1/q}$,
$\mathbb{L}_{p}(Q_T)=\mathbb{L}_{p,p}(Q_T)$,
$\|u\|_{p,p,T}=\|u\|_{p,T}$.

Let $\varrho$ be a positive function on $\mathbb{R}^d$. By
$\mathbb{L}_{p,\varrho}(\mathbb{R}^d)$
($\mathbb{L}_{p,q,\varrho}(Q_T)$) we denote the space of functions
$u$ such that $u\varrho\in \mathbb{L}_p(\mathbb{R}^d)$
($u\varrho\in\mathbb{L}_{p,q}(Q_T)$) equipped with the norm
$\|u\|_{p,\varrho}=\|u\varrho\|_p$
($\|u\|_{p,q,\varrho,T}=\|u\varrho\|_{p,q,T})$. We  write
$K\subset\subset X$ if $K$ is compact subset of $X$.
$\mathbb{L}_p^{loc}(\mathbb{R}^d)
=\bigcap_{K\subset\subset\mathbb{R}^d}\mathbb{L}_p(K)$. By
$\langle\cdot,\cdot\rangle_2$ we denote the usual inner product in
$\mathbb{L}_2(\mathbb{R}^d)$ and by
$\langle\cdot,\cdot\rangle_{2,\varrho}$ the inner product in
$\mathbb{L}_{2,\varrho}(\mathbb{R}^d)$.

$W^1_{2,\varrho}(\mathbb{R}^d)$ ($W^{0,1}_{2,\varrho}(Q_T)$) is
the Banach space consisting of all elements $u$ of
$\mathbb{L}_{2,\varrho}(\mathbb{R}^d)$
($\mathbb{L}_{2,\varrho}(Q_T)$) having generalized derivatives
$\frac{\partial u}{\partial x_i}$, $i=1,\dots,d$, in
$\mathbb{L}_{2,\varrho}(\mathbb{R}^d)$
($\mathbb{L}_{2,\varrho}(Q_T)$). If $\varrho\equiv1$ then we
denote the spaces by $W^1_{2}(\mathbb{R}^d)$ and $W^{0,1}_2(Q_T)$.
$W^{1,1}_{2,\varrho}(Q_T)$ is the subspace of
$W^{0,1}_{2,\varrho}(Q_T)$ consisting of all elements $u$ having
generalized derivatives $\frac{\partial u}{\partial t}$ in
$\mathbb{L}_{2,\varrho}(Q_T)$, $W^{1,1}_{2,0}(Q_{T})$ is the set
of all function from $W_{2}^{1,1}(Q_T)$ with compact support in
$Q_{T}$. $\mathcal{W}_\varrho=\{u\in \mathbb{L}_2
([0,T],W^1_{2,\varrho}(\mathbb{R}^d));\frac{\partial u}{\partial
t} \in\mathbb{L}_2([0,T],W^{-1}_{2,\varrho}(\mathbb{R}^d))\}$,
where $W^{-1}_{2,\varrho}(\mathbb{R}^d)$ is the dual space to
$W^{1}_{2,\varrho}(\mathbb{R}^d)$ (see
\cite{Lions,MalekNecasRokytaRuzicka} for details); if
$\varrho\equiv1$ we write $\mathcal{W}$ instead of
$\mathcal{W}_\varrho$.

By $C_{0}(Q_{T})$ $(C_{0}(\mathbb{R}^{d}))$ we denote the space of
all continuous function with compact support on $Q_{T}$
$(\mathbb{R}^{d})$ and by $C^{+}_{0}$ $(Q_{T})$
$(C^{+}_{0}(\mathbb{R}^{d}))$ the set of all positive functions
from $C_{0}(Q_{T})$ $(C_{0}(\mathbb{R}^{d}))$.

In what follows, by $C$ (or $c$) we will denote a general constant
which may vary from line to line but depends only on fixed
parameters.

\nsubsection{Preliminary results} \label{sec2}

\subsubsection{Symmetric diffusions and BSDEs}

Let $\Omega=C([0,T],\mathbb{R}^d)$ denote the space of continuous
$\mathbb{R}^d$-valued functions on $[0,T]$ equipped with the
topology of uniform convergence and let $X$ be a canonical process
on $\Omega$. It is known that given an operator $A_t$ defined by
(\ref{eq1.2}) with $a$ satisfying (\ref{eq1.1}) one can construct
a weak fundamental solution $p(s,x,t,y)$ for $A_t$ and then a
Markov family $\mathbb{X}=\{(X,P_{s,x});(s,x)\in[0,T)\times
\mathbb{R}^d\}$ for which $p$ is the transition density function,
i.e.
\[
P_{s,x}(X_t=x;0\leq t\leq s)=1,\quad
P_{s,x}(X_t\in\Gamma)=\int_{\Gamma}p(s,x,t,y)\,dy,\quad t\in(s,T]
\]
for any $\Gamma$ in a Borel $\sigma$-field $\mathcal{B}$ of
$\mathbb{R}^d$  (see \cite{A.Roz.DIFF,Stroock.DIFF}).

\begin{tw}\label{tw2.1}
For each $(s,x)\in[0,T)\times\mathbb{R}^d$, if
$[0,T)\times\mathbb{R}^d\ni (s_n,x_n)\rightarrow(s,x)$ then
$P_{s_n,x_n}\Rightarrow P_{s,x}$ weakly in $C([0,T];\mathbb{R}^d)$.
\end{tw}
\begin{dow}
Follows from the fact that $\mathbb{X}$ generates a strongly
Feller continuous Markov time-inhomogeneous semigroup on
$\mathbb{L}_{2}(\mathbb{R}^{d})$ (see \cite{A.Roz.DIFF}).
\end{dow}
\medskip

In what follows by $W$ we denote the space of all measurable
functions $\varrho:\mathbb{R}^d\rightarrow\mathbb{R}$ such that
$\varrho(x)=(1+|x|^{2})^{-\alpha}$, $x\in\mathbb{R}^d$, for some
$\alpha\ge0$.

Let $E_{s,x}$ denote expectation with respect to $P_{s,x}$.
\begin{tw}
\label{tw2.2} Let $\varrho\in W$. Then there exist $0<c\le C$
depending only on $\lambda,\Lambda$ and $\varrho$ such that
\begin{enumerate}
\item[\rm(i)]for any $\varphi\in
L_{1,\varrho}(\mathbb{R}^d)$ and $0\leq s\leq t<T$,
\begin{eqnarray*}
c\int_{\mathbb{R}^d}|\varphi(x)|\varrho(x)\,dx\leq
\int_{\mathbb{R}^d}E_{s,x}|\varphi(X_t)|\varrho(x)\,dx\leq
C\int_{\mathbb{R}^d}|\varphi(x)|\varrho(x)\,dx,
\end{eqnarray*}
\item[\rm(ii)]for any $\psi\in\mathbb{L}_{1,\varrho}(Q_T)$,
\begin{eqnarray*}
c\int_t^T\!\!\int_{\mathbb{R}^d}
|\psi(\theta,x)|\varrho(x)\,d\theta\,dx
&\le&\int_t^T\!\!\int_{\mathbb{R}^d}
E_{s,x}|\psi(\theta,X_{\theta})|\varrho(x)\,d\theta\,dx\\
&\le& C\int_t^T\!\!\int_{\mathbb{R}^d}
|\psi(\theta,x)|\varrho(x)\,d\theta\,dx,\quad t\in[s,T].
\end{eqnarray*}
\end{enumerate}
\end{tw}
\begin{dow}
Both statements follow from \cite[Proposition 5.1, Appendix]{V.Bally}, because by Aronson's
estimates there exist $0<c_1\le c_2$ depending only on
$\lambda,\Lambda$ such that
\begin{eqnarray*}
c_1\int_{\mathbb{R}^d} E|\varphi(x+X_{c_1(t-s)})|\varrho(x)\,dx
&\leq&\int_{\mathbb{R}^d}E_{s,x}|\varphi(X_t)|\varrho(x)\,dx\\
&\leq&
c_2\int_{\mathbb{R}^d}E|\varphi(x+X_{c_2(t-s)})|\varrho(x)\,dx
\end{eqnarray*}
where $E$ denotes  expectation with respect to the standard Wiener
measure on $\Omega$.
\end{dow}
\medskip

Set $\mathcal{F}^s_t=\sigma(X_u,u\in[s,t])$ and define
$\mathcal{G}$ as the completion of $\mathcal{F}^s_T$ with respect
to the family $\mathcal{P}=\{P_{s,\mu}:\mu$ is a probability
measure on $\mathcal{B}$\}, where
$P_{s,\mu}(\cdot)=\int_{\mathbb{R}^d}P_{s,x}(\cdot)\,\mu(dx)$, and
define $\mathcal{G}^s_t$ as the completion of $\mathcal{F}^s_t$ in
$\mathcal{G}$ with respect to $\mathcal{P}$.

From \cite[Theorem 2.1]{A.Roz.BSDE} it follows that there exist a
martingale additive functional locally of finite energy
$M=\{M_{s,t}:0\le s\le t\le T\}$ of $\mathbb{X}$ and a continuous
additive functional locally of zero energy $A=\{A_{s,t}:0\le s\le
t\le T\}$ of $\mathbb{X}$ such that
\begin{equation}\label{dekompozycja dyfuzji}
X_t-X_s=M_{s,t}+A_{s,t},\quad t\in[s,T],\quad P_{s,x}\mbox{-}a.s.
\end{equation}
for each $(s,x)\in[0,T)\times\mathbb{R}^d$. Moreover, the above
decomposition is unique and for each
$(s,x)\in[0,T)\times\mathbb{R}^d$ the process $M_{s,\cdot}$ is a
$(\{\mathcal{G}^s_t\}_{t\in[s,T]},P_{s,x})$-square-integrable
martingale on $[s,T]$ with the co-variation process given by
\[
\langle M^i_{s,\cdot},M^j_{s,\cdot}\rangle_t= \int_s^ta_{ij}(\theta,
X_\theta)d\theta,\quad t\in[s,T],\quad i,j=1,...,d
\]
(see \cite{A.Roz.BSDE} for details).

We now formulate  definitions of backward stochastic differential
equation (BSDE) and reflected BSDE (RBSDE) associated with
$\mathbb{X}$ and recall some known results on such equations to be
used further on.

Write
\[
B_{s,t}=\int_s^t\sigma^{-1}(\theta,X_\theta)\,dM_{s,\theta},\quad
t\in[s,T],
\]
where $M$ is the additive functional of the decomposition
(\ref{dekompozycja dyfuzji}). Notice that
$\{B_{s,t}\}_{t\in[s,T]}$ is a Wiener process.
\medskip\\
{\em Definition}  A pair $(Y^{s,x},Z^{s,x})$ of processes on $[s,T]$
is a solution of BSDE$(\varphi,f)$ (associated with $(X,P_{s,x})$)
if
\begin{enumerate}
\item[(i)] $Y^{s,x},Z^{s,x}$ are $\{\mathcal{G}^s_t\}$-adapted,
\item[(ii)] $Y^{s,x}_t=\varphi(X_T)+\int_t^T
f(\theta,X_\theta,Y^{s,x}_\theta,Z^{s,x}_\theta)\,d\theta
-\int_t^TZ^{s,x}_\theta\,dB_{s,\theta}$, $t\in[s,T]$,
$P_{s,x}$-a.s.,
\item[(iii)]$E_{s,x}\int_s^T|Z_t^{s,x}|^2\,dt<\infty,\, E_{s,x}\sup_{s\le t\le T}|Y^{s,x}_{t}|^{2}<\infty$.
\end{enumerate}
{\em Definition} A triple $(Y^{s,x},Z^{s,x},K^{s,x})$ of processes
on $[s,T]$ is a solution of RBSDE$(\varphi,f,h)$ (associated with
$(X,P_{s,x})$) if
\begin{enumerate}
\item[(i)] $Y^{s,x},Z^{s,x},K^{s,x}$ are
$\{\mathcal{G}^s_t\}$-adapted,
\item[(ii)] $Y^{s,x}_t\geq h(t,X_t)$, $t\in[s,T]$, $P_{s,x}$-a.s.,
\item[(iii)] $Y^{s,x}_t=\varphi(X_T)
+\int_t^Tf(\theta,X_\theta,Y^{s,x}_\theta,Z^{s,x}_\theta)\,d\theta
+K^{s,x}_T-K^{s,x}_t -\int_t^TZ^{s,x}_\theta\,dB_{s,\theta}$,
$t\in[s,T]$, $P_{s,x}$-a.s.,
\item[(iv)] $E_{s,x}\int_s^T|Z^{s,x}_t|^2\,dt<\infty, \, E_{s,x}\sup_{s\le t\le T}|Y^{s,x}_{t}|^{2}<\infty$,
\item[(v)] $K^{s,x}$ is a continuous increasing process such that
$K^{s,x}_s=0$, $E_{s,x}|K^{s,x}_T|^2<\infty$ and
$\int_s^T(Y^{s,x}_t-h(t,X_t))\,dK^{s,x}_t=0$, $P_{s,x}$-a.s.
\end{enumerate}

Observe that $\{\mathcal{G}^s_t\}$ need not coincide with the
natural filtration generated by the Wiener process $B_{s,\cdot}$.
Consequently, due to lack of the representation theorem for
$B_{s,\cdot}$, existence of solutions of BSDE$(\varphi,f)$ does
not follow from known results for ,,usual" BSDEs.

Existence and uniqueness of solutions of BSDE$(\varphi,f)$ for
each starting point $(s,x)\in[0,T)\times\mathbb{R}^d$ was proved
in \cite{A.Roz.BSDE} under the assumption that
$\varphi\in\mathbb{L}_2(\mathbb{R}^d)$ and $f$ satisfies (H2) with
$g\in\mathbb{L}_{p,q}(Q_T)$ for some $p,q$ such that
\begin{equation}
\label{eq2.02} p,q\in(2,\infty],   \quad \frac2{q}+\frac{d}{p}<1.
\end{equation}
(see also \cite{BallyPardoux} for existence results for
quasi-every starting point $x$ proved in the case where the
forward diffusion corresponds to symmetric divergence form
operator with time-independent coefficients but not necessarily
uniformly elliptic).

Let us recall that $u$ is said to be a weak solutions of the
Cauchy problem
\begin{equation}
\label{eq4.6} \frac{\partial u}{\partial t }+A_tu=-f_u,\quad
u(T)=\varphi
\end{equation}
(PDE$(\varphi,f)$ for short) if $u\in W^{0,1}_{2,loc}(Q_{T})\cap
C([0,T],\mathbb{L}^{loc}_{2}(\mathbb{R}^{d}))$ and for any
$\eta\in W^{1,1}_{2,0}(Q_{T})$,
\begin{eqnarray*}
&&\int_{t}^{T}\langle u(s),\frac{\partial\eta}{\partial
s}(s)\rangle_{2}\,ds+\frac12\int_{t}^{T}\langle a(s)\nabla
u(s),\nabla \eta(s)\rangle_{2}\,ds=\int_{t}^{T}\langle
f_u(s),\eta(s)\rangle_{2}\,ds\nonumber \\
&&\qquad+\langle \varphi,\eta(T)\rangle_{2}-\langle
u(t),\eta(t)\rangle_{2},\quad t\in[0,T].
\end{eqnarray*}
It is well known that if $\varphi\in\mathbb{L}_2(\mathbb{R}^d)$,
$g\in\mathbb{L}_2(Q_T)$ then there exists a unique weak solution
of PDE$(\varphi,f)$ (see, e.g. \cite{Lad}).

The next theorem strengthens slightly results proved in
\cite{A.Roz.BSDE}.

\begin{stw}\label{stw2.3}
Assume that (H1)-(H3) are satisfied with
$\varphi\in\mathbb{L}_2(\mathbb{R}^d)$, $g\in\mathbb{L}_2(Q_T)$.
\begin{enumerate}
\item[\rm(i)]If
\begin{equation}
\label{n2.2} \forall_{K\subset\subset [0,T)\times
\mathbb{R}^{d}}\quad \sup_{(s,x)\in K}
E_{s,x}\int_s^T|g(t,X_t)|^2\,dt<\infty
\end{equation}
then there exists a unique weak solution $u\in W^{0,1}_2(Q_T)\cap
C([0,T)\times\mathbb{R}^d)$ of PDE$(\varphi,f)$ and for each
$(s,x)\in[0,T)\times\mathbb{R}^d$ the pair
\begin{equation}
\label{eq02.4} (Y^{s,x}_t,Z^{s,x}_t)=(u(t,X_t),\sigma\nabla
u(t,X_t)),\quad t\in[s,T]
\end{equation}
is a unique solution of BSDE$(\varphi,f)$.
\item[\rm(ii)]There exists a version $u$ of a weak solution of
PDE$(\varphi,f)$ such that if
\begin{equation}
\label{eq2.05}
E_{s,x}\int_s^T|g(t,X_t)|^2\,dt<\infty
\end{equation}
for some $(s,x)\in[0,T)\times\mathbb{R}^d$ then the pair
(\ref{eq02.4}) is a unique solution of BSDE$(\varphi,f)$.
\end{enumerate}
\end{stw}
\begin{dow}
Let $\bar u\in W^{0,1}_2(Q_T)$ be a weak solution of the problem
(\ref{eq4.6}) and let
\[
\|\bar u\|^2_{\mathcal{W}_{2}(x,s,T)}=E_{s,x}\int^T_s(|\bar
u(t,X_t)|^2 +|\nabla\bar u(t,X_t)|^2)\,dt.
\]
From the proof of \cite[Theorem 6.1]{A.Roz.BSDE} it follows that
under (\ref{n2.2}) for every
$K\subset\subset[0,T)\times\mathbb{R}^d$,
\begin{equation}
\label{eq2.07} \sup_{(s,x)\in K}\|\bar
u\|_{\mathcal{W}_{2}(x,s,T)}<\infty.
\end{equation}
For $n,m\in\mathbb{N}$ let $u_{nm}\in W^{0,1}_2(Q_T)\cap
C([0,T)\times\mathbb{R}^d)$ be a weak solution of the Cauchy
problem
\[
\left(\frac{\partial}{\partial t} +A_t\right)u_{nm} =f_{\bar
u}^+\wedge m-f_{\bar u}^-\wedge n,\quad u_{nm}(T)=\varphi.
\]
By \cite[Proposition 5.1]{A.Roz.BSDE} the pair
$(u_{nm}(t,X_t),\sigma\nabla u_{nm}(t,X_t))$, $t\in[s,T]$, is a
solution of BSDE$(\varphi, f^+_{\bar u}\wedge m-f^-_{\bar u}\wedge
n)$. Using It\^{o}'s formula and performing standard calculations
we conclude that there is $C>0$ not depending on $n,m$ such that
\begin{eqnarray}\label{eq2.6}
&&E_{s,x}\sup_{s\le t\le T}
|u_{nm}(t,X_t)|^2+E_{s,x}\int_s^T|\sigma\nabla
u_{nm}(t,X_t)|^2\,dt\nonumber\\
&&\qquad \le C\left(E_{s,x}|\varphi(X_T)|^2
+E_{s,x}\int_s^T|g(t,X_t)|^2\,dt+\|\bar
u\|^2_{\mathcal{W}_{2}(x,s,T)} \right).
\end{eqnarray}
From comparison results (see \cite[Theorem
4.1.4]{BensoussanLions}) and the fact that $u_{nm}$ are continuous
it follows that for any fixed $n$ the sequence
$\{u_{nm}\}_{m\in\mathbb{N}}$ is increasing. Hence, for each
$n\in\mathbb{N}$ there is $u_n$ such that $u_{nm}\uparrow u_n$ as
$m\rightarrow\infty$. Moreover, by well known convergence theorems
(see \cite[Theorem 3.4.5]{Lad}), $u_{nm}\rightarrow u_n$ in
$W^{0,1}_2(Q_T)$ and $u_n$ is a weak solution of the problem
\[
\left(\frac{\partial}{\partial t}+A_t\right)u_n=f^+_{\bar
u}-f^-_{\bar u}\wedge n,\quad u_n(T)=\varphi.
\]
If (\ref{n2.2}) is satisfied, then from (\ref{eq2.07}),
(\ref{eq2.6}) and Nash's continuity theorem (see \cite{Aro}) it
follows that $\{u_{nm}\}_{m\in\mathbb{N}}$ is equicontinous in
every compact subset of $[0,T)\times\mathbb{R}^d$. Therefore the
functions $u_n$ are continuous on $[s,T)\times\mathbb{R}^d$. Using
once again It\^{o}'s formula we deduce that for any
$k,l,n\in\mathbb{N}$,
\begin{eqnarray}
\label{eq2.7} &&E_{s,x}|(u_{nk}-u_{nl})(t,X_t)|^2 +E_{s,x}\int_s^T
|\sigma\nabla(u_{nk}-u_{nl})(t,X_t)|^2\,dt\nonumber\\
&&\qquad\le C\left(E_{s,x}\int_s^T |(f^+_{\bar u}\wedge
k-f^+_{\bar u}\wedge l)
(t,X_t)|^2\,dt\right)^{1/2}\nonumber \\
&&\qquad\qquad\qquad\qquad\times\left(E_{s,x}\int_s^T
|(u_{nk}-u_{nl})(t,X_t)|^2\,dt\right)^{1/2}
\end{eqnarray}
for all $t\in[s,T]$. By (H2) and (\ref{eq2.07}), (\ref{eq2.6}) the
first term on the right-hand side of (\ref{eq2.7}) is bounded
uniformly in $k,l$. Due to (\ref{eq2.07}), (\ref{eq2.6}) and the
estimate $|u_{nk}|\le|u_{n1}|+|u_{n}|$ we  may apply the Lebesgue
dominated convergence theorem to conclude that the second term
converges to zero as $k,l\rightarrow0$. By the above,
\[
E_{s,x}|(u_{nm}-u_n)(t,X_t)|^2 +E_{s,x}\int_s^T
|\sigma\nabla(u_{nm}-u_n)(t,X_t)|^2\,dt\rightarrow0
\]
as $m\rightarrow\infty$. Using this it is easy to see that the
pair $(u_n(t,X_t),\sigma\nabla u_n(t,X_t))$, $t\in[s,T]$ is a
solution of BSDE$(\varphi,f^+_{\bar u}-f^-_{\bar u}\wedge n)$.
Therefore, (\ref{eq2.6}) holds for $u_{nm}$ replaced by $u_n$ and
(\ref{eq2.7}) holds for $u_{nk},u_{nl}$ replaced by $u_k,u_l$ and
$f^+_{\bar u}$ replaced by $f^-_{\bar u}$. Using once again
(\ref{eq2.07}) and Nash's continuity theorem we conclude that
$u_{n}$ is equicontinuous in every compact subset of $[0,T)\times
\mathbb{R}^{d}$. Therefore, by comparison results, ${u_{n}}$ is
decreasing and there is $u\in C([0,T)\times\mathbb{R}^d)$ such
that $u_n\downarrow u$. Since $f^+_{\bar u}-f^-_{\bar u}\wedge
n\rightarrow f_{\bar u}$ in $\mathbb{L}_2(Q_T)$, it follows that
$u$ is a weak solution of the Cauchy problem
$(\frac{\partial}{\partial t}+A_t)u=f_{\bar u}$, $u(T)=\varphi$.
By uniqueness, $u$ is a version of $\bar u$. Finally, using the
mentioned above analogues of (\ref{eq2.6}), (\ref{eq2.7}) we prove
in much the same way as above that the pair (\ref{eq02.4}) is a
solution of BSDE$(\varphi,f)$, which completes the proof of (i).

To prove (ii), we first observe that using continuity of
${u_{nm}}$ and the fact that $\{u_{nm}\}$ is decreasing for every
fixed $n$ and increasing for every fixed $m$ we can still show
that $\{u_{n}\}$ is decreasing. Therefore $\{u_n\}$ converges
pointwise to some version $u$ of $\bar u$. If (\ref{eq2.05}) is
satisfied for some $(s,x)\in [0,T)\times \mathbb{R}^{d}$ then
$\|\bar u\|_{\mathcal{W}_{2}(x,s,T)}<\infty$. Therefore we can
use (\ref{eq2.6}), (\ref{eq2.7}) to conclude as before that
$(u(t,X_t),\sigma\nabla u(t,X_t))$, $t\in[s,T]$, is a solution of
BSDE$(\varphi,f)$ associated with $(X,P_{s,x})$.
\end{dow}

\begin{tw}
\label{tw2.4} Assume that (H1)--(H3) are satisfied with
$\varphi\in \mathbb{L}_{2}(\mathbb{R}^{d})$,
$h,g\in\mathbb{L}_2(Q_{T})$ and
\begin{equation}
\label{eq3.1} E_{s,x}\sup_{s\leq t\leq T}
|h^{+}(t,X_t)|^2+E_{s,x}\int_s^T|g(t,X_t)|^2\,dt<\infty
\end{equation}
for some $(s,x)\in [0,T)\times\mathbb{R}^d$. Then the
RBSDE$(\varphi,f,h)$ associated with $(X,P_{s,x})$ has a unique
solution $(Y^{s,x},Z^{s,x},K^{s,x})$. Moreover, if the pair
$(Y^{s,x,n}_t,Z^{s,x,n})$, $n\in\mathbb{N}$, is a solution of
BSDE$(\varphi,f+n(y-h)^{+})$, then
\begin{eqnarray}
\label{zbieznosci} &&E_{s,x}\sup_{s\leq t\leq
T}|Y^{s,x,n}_t-Y^{s,x}_t|^2
+E_{s,x}\int_s^T|Z^{s,x,n}_t-Z^{s,x}_t|^2\,dt\nonumber\\
&&\qquad+ E_{s,x}\sup_{s\leq t\leq T}
|K^{s,x,n}_t-K^{s,x}_t|^2\rightarrow0,
\end{eqnarray}
where
\[
K_t^{s,x,n}=\int_s^tn(Y^{s,x,n}_\theta-h(\theta,X_\theta))d\theta,\quad
t\in[s,T],\quad P_{s,x}\mbox{-}a.s.
\]
Finally, there is $C>0$ depending neither on $n,m\in\mathbb{N}$
nor on $s,x$ such that
\begin{eqnarray}
\label{eq2.4} &&E_{s,x}\sup_{s\leq t\leq T}
|Y^{s,x,n}_t|^2+E_{s,x}\int_s^T|Z^{s,x,n}_t|^2\,dt
+E_{s,x}|K^{s,x,n}_T|^2 \nonumber\\
&&\quad\leq C\left(E_{s,x}|\varphi(X_T)|^2 +E_{s,x}\sup_{s\leq
t\leq T} |h^+(t,X_t)|^2 +E_{s,x}\int_s^T|g(t,X_t)|^2\,dt\right)
\end{eqnarray}
and
\begin{eqnarray}\label{ciaglosc}
&&E_{s,x}\sup_{s\leq t\leq
T-2\delta}|Y^{s,x,n}_t-Y^{s,x,m}_t|^2\nonumber\\
&&\quad\le
C\left(E_{s,x}|Y^{s,x,n}_{T-\delta}-Y^{s,x,m}_{T-\delta}|^2
+E_{s,x}\int_s^{T-\delta}
(Y^{s,x,n}_t-h(t,X_t))^-\,dK^{s,x,m}_t\right.\nonumber\\
&&\left.\qquad\qquad\qquad
+E_{s,x}\int_s^{T-\delta}(Y^{s,x,m}_t-h(t,X_t))^-\,dK^{s,x,n}_t\right)
\end{eqnarray}
for every $\delta\in[0,T-s]$.
\end{tw}
\begin{dow}
From Proposition \ref{stw2.3} we know that for each
$n\in\mathbb{N}$ there exists a unique solution of
BSDE$(\varphi,f+n(y-h)^{+})$. To prove
(\ref{zbieznosci})--(\ref{ciaglosc}) it suffices to repeat step by
step arguments from the proofs of corresponding results in
\cite{EKPPQ}.
\end{dow}
\medskip

Let us remark that both terms in (\ref{eq3.1}) are bounded
uniformly in $(s,x)\in K$ for every
$K\subset\subset[0,T)\times\mathbb{R}^d$ if $h,g$ satisfy the
polynomial growth condition or $h$ satisfies the polynomial growth
condition and $g\in\mathbb{L}_{p,q,\varrho}(Q_{T})$ with $p,q$
satisfying (\ref{eq2.02}) and $\varrho\in W$. The first statement
is an immediate consequence of Proposition \ref{doob} proved in
Section \ref{sec3}. Sufficiency of the second condition on $g$
follows from H\"older's inequality and upper Aronson's estimate on
the transition density $p$ (see \cite{Aro}).

Observe also that if $g\in\mathbb{L}_{2,\varrho}(Q_T)$ then
(\ref{eq2.05}) holds for a.e. $(s,x)\in[0,T)\times\mathbb{R}^d$
because by Theorem \ref{tw2.2},
\[
\int_{0}^{T}\left(\int_{\mathbb{R}^{d}}
(E_{s,x}\int_{s}^{T}|g(t,X_{t})|^2\,dt\varrho^{2}(x))\,dx\right)ds\le
C\|g\|_{2,\varrho}^{2}
\]

\begin{lm}\label{oszacowanie przyrostu}
If $(Y^{s,x,i}_t,Z^{s,x,i}_t,K^{s,x,i}_t)$, $i=1,2$,  is  a solution
of RBSDE$(\xi,f,h^{i})$ then for every $\delta\in[0,T-s]$,
\begin{eqnarray*}
&&E_{s,x}\sup_{s\leq t\leq T-\delta}|Y^{s,x,1}_t-Y^{s,x,2}_t|^2
+E_{s,x}\int_s^{T-\delta}|Z^{s,x,1}_t-Z^{s,x,2}_t|^2\,dt\\
&&\qquad\qquad+E_{s,x}\sup_{s\leq t\leq T-\delta}|K^{s,x,1}_t-K^{s,x,2}_t|^2\\
&&\qquad\leq C\left(E_{s,x}\sup_{s\leq t\leq T-\delta}
|h^{1}(t,X_t)-h^2(t,X_t)|^2
+E_{s,x}|Y^{s,x,1}_{T-\delta}-Y^{s,x,2}_{T-\delta}|^2\right).
\end{eqnarray*}
\end{lm}
\begin{dow}
See \cite{EKPPQ}.
\end{dow}

\subsubsection{Obstacle problem}
\label{ssec2.2}

In this subsection we formulate precisely our  definition of
solutions of the obstacle problem and compare it to the well known
definitions of solutions in the sense of variational inequalities.
We prove also a priori estimates for solutions and some additional
technical results which will be needed in the next section.

In the paper we will use the following notion of the capacity of
$E\subset\subset\breve{Q}_T$:
\[
cap_{\breve{Q}_{T}}(E)=\inf\{\int_{\breve{Q}_{T}}
(|\frac{\partial\eta}{\partial
t}(t,x)|^2+|\nabla\eta(t,x)|^{2})\,dt\,dx:\eta\in
C_{0}^{\infty}(\breve{Q}_{T}), \eta\geq\mathbf{1}_{E}\}.
\]
In the standard way we can extend the above capacity to external
capacity for arbitrary subset $E\subset\breve{Q}_{T}$. It is known
that $cap_{\breve{Q}_{T}}$ is  the Choquet capacity (see Chapter 2
in  \cite{Fuk.O.Tak.}).

In the remainder of the paper the abbreviation ``q.e."  means
``except for a set of capacity zero".

Throughout the subsection we assume that $\varrho\in W$ and
(H1)--(H3) are satisfied.
\medskip\\
{\em Definition} We say that a pair $(u,\mu)$, where $\mu$ is a
Radon  measure on $Q_T$ and $u:Q_T\rightarrow\mathbb{R}$ is a
measurable function defined up to the sets of $\mu$-measure zero,
is a weak solution of the obstacle problem (\ref{eq1.3}) with data
$\varphi,f,h$ (OP$(\varphi,f,h)$ for short) if
\begin{enumerate}
\item[(a)]$u\in
W^{0,1}_{2,loc}(Q_{T})\cap
C([0,T],\mathbb{L}^{loc}_{2}(\mathbb{R}^{d}))$ and for any $\eta\in
W^{1,1}_{2,0}(Q_{T})$,
\begin{eqnarray}
\label{eq2.1} &&\int_{t}^{T}\langle
u(s),\frac{\partial\eta}{\partial
s}(s)\rangle_{2}\,ds+\frac12\int_{t}^{T}\langle a(s)\nabla
u(s),\nabla \eta(s)\rangle_{2}\,ds=\int_{t}^{T}\langle
f_u(s),\eta(s)\rangle_{2}\,ds\nonumber \\
&&\qquad+\int_{t}^{T}\!\!\int_{\mathbb{R}^{d}}\eta \,d\mu+\langle
\varphi,\eta(T)\rangle_{2}-\langle u(t),\eta(t)\rangle_{2},\quad
t\in[0,T],
\end{eqnarray}
\item[(b)] $u\geq h$ on $Q_T$,
\item[(c)] $\int_{Q_T}(u-h)\xi\,d\mu=0$ for all $\xi\in
C^{+}_{0}(Q_{T})$,
\item[(d)] $\mu(\{t\}\times\mathbb{R}^d)=0$ for every $t\in[0,T]$.
\end{enumerate}

Some comments on the above definition are in order. In the next
lemma we will show that (a) forces $\mu_{|\breve{Q}_{T}}\ll
cap_{\breve{Q}_{T}}$, which together with (d) and the well known
fact that elements of $W_{2}^{1,1}(Q_{T})$ are defined up to
subsets of $\breve Q_T$ of zero capacity (see, e.g.,
\cite{B.M.W.C.,Evance,M}) ensures that the integral
$\int_{t}^{T}\!\!\int_{\mathbb{R}^{d}}\eta\,d\mu$ is correctly
defined. We shall see that (a) implies that
$\mu(\{s\}\times\mathbb{R}^{d})=0$ for $s\in(0,T)$, so instead of
(d) we could impose the condition
$\mu(\{0,T\}\times\mathbb{R}^{d})=0$. The condition
$\mu(\{T\}\times\mathbb{R}^{d})=0$ is also necessary for the
terminal condition $u(T)=\varphi(T)$ to hold, and a fortiori, for
uniqueness of the solution of the obstacle problem. Notice also
that the integral in condition (c) is well defined because
$u-h\ge0$.

Let us remark that our definition of the obstacle problem is
similar to that in stochastic case (condition (d) may be viewed as
an analytical counterpart to continuity of the process $K^{s,x}$).
Notice also that if the obstacle $h$ is constant, then the above
definition coincides with the one adopted in \cite{Mokrane} (in
\cite{Mokrane} exclusively constant obstacles are considered; this
implies that $\mu$ is absolutely continuous with respect to the
Lebesgue measure, so no problems arises with the definition of an
obstacle problem).

\begin{lm}
\label{lemat o rownowaznosci} If $u\in W^{0,1}_{2,loc}(Q_{T})\cap
C([0,T],\mathbb{L}^{loc}_{2}(\mathbb{R}^{d}))$ and  the pair
$(u,\mu)$ satisfies
\begin{eqnarray}
\label{eq2.2} &&\int_{0}^{T}\langle
u(t),\frac{\partial\eta}{\partial t}(t)\rangle_{2}\,dt
+\frac12\int_{0}^{T}\langle
a(t)\nabla u(t),\nabla\eta(t)\rangle_{2}\,dt \nonumber\\
&&\qquad=\int_{0}^{T}\langle f_u(t),\eta(t)\rangle_{2}\,dt
+\int_{Q_T}\eta\,d\mu +\langle\varphi,\eta(T)\rangle_{2}
\end{eqnarray}
for every $\eta\in W^{1,1}_{2}(Q_{T})\cap C_0(Q_T)$ such that
$\eta(0)\equiv0$, then

\begin{enumerate}
\item[\rm(i)]$\mu_{|\breve{Q}_{T}}\ll cap_{\breve{Q}_{T}}$\,,
\item[\rm(ii)]$\mu(\{t\}\times\mathbb{R}^{d})=0$ for every $t\in(0,T)$,
\item[\rm(iii)] $u(T)=\varphi$ if and only if
$\mu(\{T\}\times\mathbb{R}^{d})=0$,
\item[\rm(iv)] if $\mu(\{0,T\}\times \mathbb{R}^{d})=0$ then
(\ref{eq2.1}) holds for $\eta\in W^{1,1}_{2,0}(Q_{T})$.
\end{enumerate}
\end{lm}
\begin{dow}
Fix $E\subset\subset\breve{Q}_{T}$ and choose positive $\eta\in
C_{0}^{\infty}(\breve{Q}_{T})$ such that $\eta\geq
\mathbf{1}_{E}$. Then, by (\ref{eq2.2}),
\[
\mu(E)\le \int_{0}^{T}\langle u(t),\frac{\partial\eta}{\partial
t}(t)\rangle_{2}\,dt+\frac12\int_{0}^{T}\langle a(t)\nabla
u(t),\nabla \eta(t)\rangle_{2}\,dt-\int_{0}^{T}\langle
f_u(t),\eta(t)\rangle_{2}\,dt,
\]
and hence, by Gagliardo-Nirenberg-Sobolev inequality,
\[
\mu(E)\le C (cap_{\breve{Q}_{T}}(E))^{1/2} (\|u\|_{2,T}+\|\nabla
u\|_{2,T}+\|f_{u}\|_{2,T}),
\]
which shows (i). Now, fix $s\in(0,T)$ and consider the sequence of
functions $\{\eta^{n,s}\}$ defined  by
\[
\eta^{n,s}(t,x)
=\left\{\begin{array}{ll} 0,& t\in[0,s_{n}], \\
\frac{\eta(s,x)}{s-s_{n}}(t-s_{n}),& t\in(s_{n},s),\\
\eta(t,x),& t\in[s,T],
\end{array}\right.
\]
where $\{s_n\}\subset(0,s)$ is a sequence such that $s_n\uparrow s$.
Observe that $\eta^{n,s}\rightarrow
\mathbf{1}_{[s,T]\times\mathbb{R}^{d}}\eta$, $\nabla
\eta^{n,s}\rightarrow\mathbf{1}_{[s,T]\times\mathbb{R}^{d}}\nabla\eta$
and
\[
\frac{\partial\eta}{\partial t}^{n,s}(t,x)
=\left\{\begin{array}{ll} 0,& t\in[0,s_{n}],\\
\frac{\eta(s,x)}{s-s_{n}},& t\in(s_n,s),\\
\frac{\partial\eta}{\partial t}(t,x),& t\in[s,T].
\end{array}\right.
\]
From (\ref{eq2.2}) with $\eta$ replaced by $\eta^{n,s}$ we have
\begin{eqnarray*}
&&\frac{1}{s-s_n}\int_{s_n}^{s}\langle
u(t),\eta(t)\rangle_{2}\,dt+\int_{s}^{T}\langle u(t),\frac{\partial
\eta}{\partial t}(t)\rangle_{2}\,dt+\frac12\int_{0}^{T}\langle
a(t)\nabla u(t),\nabla\eta^{n,s}(t)\rangle\,dt\\
&&\qquad =\int_{0}^{T}\langle
f_u(t),\eta^{n,s}(t)\rangle_{2}\,dt+\int_{Q_T}\eta^{n,s}\,d\mu
+\langle\varphi,\eta^{n,s}(T)\rangle_{2}.
\end{eqnarray*}
Letting $n\rightarrow\infty$ and using the fact that $u\in
C([0,T],\mathbb{L}^{loc}_{2}(\mathbb{R}^{d}))$ we get
(\ref{eq2.1}) for every $\eta\in W^{1,1}_{2}(Q_{T})\cap C_0(Q_T),
t\in(0,T]$. In particular,  for any positive $\eta\in
W^{1,1}_{2}(Q_{T})\cap C_0(Q_T)$ and any $0<h<s\le T$ we have
\begin{eqnarray*}
\int_{s-h}^{s}\!\int_{\mathbb{R}^{d}}\eta\,d\mu
&=&\int_{s-h}^{s}\langle u(t),\frac{\partial\eta}{\partial
t}(t)\rangle_{2}\,dt+\frac12\int_{s-h}^{s}\langle
a(t)\nabla u(t),\nabla \eta(t)\rangle_{2}\,dt\\
&&-\int_{s-h}^{s}\langle f_u(t),\eta(t)\rangle_{2}\,dt-\langle
u(s),\eta(s)\rangle_{2} +\langle u(s-h),\eta(s-h)\rangle_{2},
\end{eqnarray*}
so letting $h\downarrow0$ and using continuity of $t\mapsto u(t)$
in $\mathbb{L}_2(\mathbb{R}^{d})$ we get (ii) and (iii). To show
(iv) we assume that $\eta\in W^{1,1}_{2,0}(Q_{T})$ and consider a
sequence $\{\eta_{n}\}\subset W^{1,1}_{2}(Q_{T})\cap C_0(Q_T)$
such that $\eta_{n}\rightarrow\eta$ in $ W^{1,1}_{2}(Q_{T})$ and
quasi-everywhere in $\breve{Q}_{T}$. From (i) and the assumption
in (iv) it follows that $\{\eta_n\}$ converges $\mu$-a.e. in
$Q_{T}$ as well. From (\ref{eq2.1}) applied to $|\eta_n-\eta_m|$
we conclude  that $\{\eta_{n}\}$ is a Cauchy sequence in
$\mathbb{L}_1([t,T]\times \mathbb{R}^{d},\mu)$ for every $t\in
(0,T]$. Therefore (\ref{eq2.1}) is satisfied for any $\eta\in
W^{1,1}_{2,0}(Q_{T})$ and $t\in (0,T]$. Clearly, if
$\mu(\{0\}\times\mathbb{R}^d)=0$, then it is satisfied also for
$t=0$.
\end{dow}
\medskip

In what follows, given some function $u:Q_T\rightarrow\mathbb{R}^{d}$ we
will extend it in a natural way to the function on
$[-T,2T]\times\mathbb{R}^{d}$, still denoted by $u$, by putting
\begin{displaymath}
u(t,x)=\left\{\begin{array}{ll} u(-t,x),& t\in[-T,0],\\
u(t,x),& t\in[0,T],\\
u(2T-t,x),& t\in[T,2T].
\end{array}\right.
\end{displaymath}
For $\varepsilon>0$ set
\[
u_{\varepsilon}(t,x)=\frac{1}{\varepsilon}\int_{0}^{\varepsilon}u(t-s,x)\,ds,
\quad (t,x)\in[0,T]\times\mathbb{R}^{d}
\]
and note that if $u\in
C([0,T],\mathbb{L}^{loc}_{2}(\mathbb{R}^{d}))\cap
W^{0,1}_{2,loc}(Q_{T})$ then $u_{\varepsilon}\in
W^{1,1}_{2,loc}(Q_{T})$, $t\mapsto
u_{\varepsilon}(t)\in\mathbb{L}^{loc}_{2}(\mathbb{R}^{d})$ is
differentiable,  $\nabla u_{\varepsilon}\rightarrow \nabla u$ in
$\mathbb{L}^{loc}_{2}(Q_{T})$ and $u_{\varepsilon}(t)\rightarrow
u(t)$ in $\mathbb{L}^{loc}_{2}(\mathbb{R}^{d})$ for every
$t\in[0,T]$.

\begin{lm}\label{lm2.8}
If $(u,\mu)$ satisfies (a),(d), then for any $\eta\in
W^{1,1}_{2,0}(Q_{T})$ and $t\in(0,T)$,
\begin{eqnarray}
\label{eq2.3} &&\int_{t}^{T}\langle
u_{\varepsilon}(s),\frac{\partial \eta}{\partial
s}(s)\rangle_{2}\,ds +\frac12\int_{t}^{T}\langle a(s)\nabla
u_{\varepsilon}(s),\nabla\eta(s)\rangle_2\,ds
\nonumber\\
&&\quad=\int_{t}^{T}\langle
f_{u,\varepsilon}(s),\eta(s)\rangle_{2}\,ds
+\frac{1}{\varepsilon}\int_{0}^{\varepsilon}\left(\int_{t-\theta}^{T-\theta}
\!\!\int_{\mathbb{R}^{d}}\eta(s+\theta,x)\,d\mu(s,x)\right)d\theta\nonumber\\
&&\qquad+\langle u_{\varepsilon}(T),\eta(T)\rangle_{2} -\langle
u_{\varepsilon}(t),\eta(t)\rangle_2
\end{eqnarray}
for all sufficiently small $\varepsilon>0$.
\end{lm}
\begin{dow}
Using Fubini's theorem and (\ref{eq2.1}) we obtain
\begin{eqnarray*}
&&\int_t^T\langle u_{\varepsilon}(s),\frac{\partial\eta}{\partial
s}(s)\rangle_{2}\,ds
=\frac{1}{\varepsilon}\int_{0}^{\varepsilon}\left(\int_{t-\theta}^{T-\theta}
\langle u(s),\frac{\partial\eta}{\partial s}
(s+\theta)\rangle_2\,ds\right)\,d\theta\\
&&\quad=-\frac{1}{2\varepsilon}\int_{0}^{\varepsilon}
\left(\int_{t-\theta}^{T-\theta}\langle a(s)\nabla
u(s),\nabla\eta(s+\theta)\rangle_2\,ds\right)\,d\theta\\
&&\qquad+\frac{1}{\varepsilon}\int_{0}^{\varepsilon}
\left(\int_{t-\theta}^{T-\theta} \langle
f_u(s),\eta(s+\theta)\rangle_2\,ds
+\int_{t-\theta}^{T-\theta}\!\!\int_{\mathbb{R}^{d}}
\eta(s+\theta,x)\,d\mu(s,x)\right)\,d\theta\\
&&\qquad+\frac{1}{\varepsilon}\int_{0}^{\varepsilon}(\langle
u(T-\theta),\eta(T)\rangle_2 -\langle
u(t-\theta),\eta(t)\rangle_2)\,d\theta,
\end{eqnarray*}
from which (\ref{eq2.3}) follows.
\end{dow}

\begin{stw}\label{stw2.8}
If  $(u,\mu)$ satisfies (a), (d) and $u\in C(\breve Q_T)$ then
$\int_{Q_{T}}\xi^2|u|\,d\mu<\infty$ for any $\xi\in C_0^{1}(Q_T)$.
Moreover,
\begin{eqnarray}
\label{eq2.06} &&\|u(t)\xi\|^{2}_{2} +\int_{t}^{T}\langle
a(s)\nabla u(s),\nabla(u\xi^2)(s)\rangle_{2}\,ds\nonumber\\
&&\qquad=\|\varphi\xi\|^{2}_{2} +2\int_{t}^{T}\langle
f_u(s),u(s)\xi^2\rangle_{2}\,ds
+2\int_{t}^{T}\!\!\int_{\mathbb{R}^{d}}\xi^2u\,d\mu
\end{eqnarray}
for all $t\in[0,T]$.
\end{stw}
\begin{dow}
Let $\tau\in(0,T)$. Write
$u_{\varepsilon}^{+}=(u_{\varepsilon})^{+}$. By (\ref{eq2.3}) with
$\eta=\xi^{2}u^+_{\varepsilon}$ we have
\begin{eqnarray}
\label{eq2.5} &&\int_{t}^{\tau}\langle
u_{\varepsilon}(s),\xi^2\frac{\partial
u^{+}_{\varepsilon}}{\partial s}(s)\rangle_2\,ds
+\frac12\int_{t}^{\tau}\langle a(s)\nabla
u_{\varepsilon}(s),\nabla
(\xi^{2}u^{+}_{\varepsilon})(s)\rangle_2\,ds \nonumber\\
&&\quad=\int_{t}^{\tau}\langle
f_{u,\varepsilon}(s),\xi^2u^{+}_{\varepsilon}(s)\rangle_{2}\,ds+\langle
u_{\varepsilon}(\tau),\xi^2u_{\varepsilon}^{+}(\tau)\rangle_{2}
-\langle u_{\varepsilon}(t),\xi^2u_{\varepsilon}^{+}(t)\rangle_{2}\nonumber\\
&&\qquad+\frac{1}{\varepsilon}\int_{0}^{\varepsilon}
\left(\int_{t-s_{1}}^{\tau-s_{1}}\!\!\int_{\mathbb{R}^{d}}
\xi^2u^{+}_{\varepsilon}(s+s_{1},x)\,d\mu(s,x)\right)ds_{1}\nonumber\\
&&\quad=\int_{t}^{\tau}\langle
f_{u,\varepsilon}(s),\xi^2u^{+}_{\varepsilon}(s)\rangle_{2}\,ds
+\|\xi u_{\varepsilon}^{+}(\tau)\|_{2}^{2} -\|\xi
u_{\varepsilon}^{+}(t)\|_{2}^{2}
+\int_{Q_{T}}g^{\varepsilon}_{\xi}\,d\mu,
\end{eqnarray}
where
\begin{eqnarray*}
g^{\varepsilon}_{\xi}(s,x)
&=&\frac{1}{\varepsilon^{2}}\int_{0}^{\varepsilon}
\mathbf{1}_{[t-s_{1},\tau-s_{1}]}(s)\xi^2
\left(\int_{0}^{\varepsilon}u(s+s_{1}-s_{2},x)ds_{2}\right)^{+}ds_{1}
\end{eqnarray*}
for $s\in[t,\tau)$ and $g^{\varepsilon}_{\xi}(\tau,x)=0$. Observe
that for every $(s,x)\in[t,\tau)\times\mathbb{R}^{d}$,
\[
g^{\varepsilon}_{\xi}(s,x)=\frac{1}{\varepsilon^{2}}\int_{0}^{\varepsilon}
\xi^2\left(\int_{0}^{\varepsilon}u(s+s_{1}-s_{2},x)ds_{2}\right)^{+}ds_{1}
\]
for  sufficiently small $\varepsilon>0$. Since $|a^{+}-b^{+}|\leq
|a-b|$ for every $a,b\in\mathbb{R}$, we have
\begin{eqnarray*}
&&\left|\frac{1}{\varepsilon^{2}}\int_{0}^{\varepsilon}
\xi^2\left(\int_{0}^{\varepsilon}
u(s+s_{1}-s_{2},x)\,ds_{1}\right)^{+}\,ds_{2}
-\frac{1}{\varepsilon^{2}}\int_{0}^{\varepsilon}
\xi^2\left(\int_{0}^{\varepsilon}u(s,x)\,ds_{1}\right)^{+}\,ds_{2}\right|\\
&&\quad\leq\frac{1}{\varepsilon^{2}}\int_{0}^{\varepsilon}
\xi^{2}\left|\left(\int_{0}^{\varepsilon}
u(s+s_{1}-s_{2},x)\,ds_{1}\right)^{+}-
\left(\int_{0}^{\varepsilon}u(s,x)\,ds_{1}\right)^{+}\right|\,ds_{2}\\
&&\quad\leq\frac{1}{\varepsilon^{2}}\int_{0}^{\varepsilon}
\xi^{2}\left|\int_{0}^{\varepsilon}
u(s+s_{1}-s_{2},x)\,ds_{1}
-\int_{0}^{\varepsilon}u(s,x)\,ds_{1}\right|\,ds_{2}\\
&&\quad\leq\frac{1}{\varepsilon^{2}}
\int_0^{\varepsilon}\!\int^{\varepsilon}_0
\xi^2|u(s+s_{1}-s_{2},x)-u(s,x)|\,ds_{1}\,ds_{2},
\end{eqnarray*}
and consequently, $g^{\varepsilon}_{\xi}(s,x)\rightarrow
\xi^2u^{+}(s,x)$ for every $(s,x)\in[t,\tau)\times\mathbb{R}^{d}$
as $\varepsilon\rightarrow0$. Therefore from (\ref{eq2.5}) we
obtain
\begin{eqnarray}
\label{eq2.15} &&\frac{1}{2}\|\xi\varphi^{+}\|^{2}_{2}
-\frac{1}{2}\|\xi u^{+}(t)\|^{2}_{2}+\frac12\int_{t}^{\tau}
\langle a(s)\nabla u^{+}(s),\nabla(\xi^2u^{+})(s)\rangle_2\nonumber\\
&&\qquad=\int_{t}^{\tau}\langle f_u(s),\xi^2u^+(s)\rangle_{2}ds
+\|\xi\varphi^{+}\|^{2}_{2}-\|\xi u^{+}(t)\|^{2}_{2}\nonumber\\
&&\qquad\quad+\liminf_{\varepsilon\rightarrow0}
\int^\tau_t\!\!\int_{\mathbb{R}^d}g_{\xi,\varepsilon}(s,x)\,d\mu.
\end{eqnarray}
Hence, by Fatou's lemma,
\begin{eqnarray*}
\int^\tau_t\!\!\int_{\mathbb{R}^d}\xi^2u^{+}\,d\mu &\leq&
\frac{1}{2}\|\xi u(t)\|^{2}_{2} +\frac12\int_{t}^{\tau}|\langle
a(s)\nabla u(s),\nabla(\xi^2u)(s)\rangle_2|\,ds\\
&&+\int_{t}^{\tau}|\langle
f_u(s),\xi^2u(s)\rangle_{2}|\,ds+\frac{1}{2}\|
\xi\varphi\|^{2}_{2}.
\end{eqnarray*}
Letting  $t\downarrow0$ and $\tau\uparrow T$ we see from the above
that $\int_{Q_{T}}\xi^2u^+\,d\mu<\infty$. Analogously, putting
$\eta=\xi^2u^-_{\varepsilon}$ we show that
$\int_{Q_{T}}\xi^2u^-\,d\mu<\infty$, which completes the proof of
the first part of the lemma.  Since $|g_{\xi,\varepsilon}(s,x)|\le
C\xi^2$ on $[t,\tau]\times\mathbb{R}^d$ for some $C>0$, using the
Lebesgue dominated convergence theorem we conclude from
(\ref{eq2.15}) that (\ref{eq2.06}) is satisfied with $T$ replaced
by $\tau$ and $t\in(0,\tau]$. Because we know already that
$\int_{Q_{T}}\xi^2|u|\,d\mu<\infty$, letting $\tau\uparrow T$ and
$t\downarrow 0$ we complete the proof.
\end{dow}
\medskip

We now are ready to prove useful a priori estimates for solutions
of an obstacle problem.
\begin{stw}
\label{stw2.30} Assume (H1)--(H3) with
$\varphi\in\mathbb{L}_{2,\varrho}(\mathbb{R}^d)$,
$g\in\mathbb{L}_{2,\varrho}(Q_T)$. If $(u,\mu)$ satisfies (a) and
(d), $u\in C(\breve{Q}_{T})$, $\nabla u \in
\mathbb{L}_{2,\varrho}(Q_{T})$ and there is
$p:Q_{T}\rightarrow\mathbb{R}$ such that $p^{+}\in
W^{1,1}_{2,\varrho}(Q_{T})$ and $\int_{Q_T}(u-p)\xi d\mu\leq0$ for
all $\xi \in C^{+}_{0}(Q_{T})$ then there is $C>0$ such that
\begin{eqnarray}
\label{nierownosc1trzyStrony}
&&\sup_{t\in[0,T]}\|u(t)\|^{2}_{2,\varrho} +\int_0^T\|\nabla
u(s)\|^{2}_{2,\varrho}\,ds +\int_{Q_{T}}|u|\varrho^{2}\,d\mu +
\|\mu\|_{(W_{2,\varrho}^{1,1}(Q_{T}))^{*}}\nonumber\\
&&\qquad\leq C\left(\|\varphi\|^{2}_{2,\varrho}
+\sup_{t\in[0,T]}\|p^{+}(t)\|^{2}_{2,\varrho}\right.\nonumber\\
&&\qquad\qquad\quad\left. +\int^T_0(\|\frac{\partial
p^{+}}{\partial s}(s)\|^{2}_{2,\varrho}+\|\nabla
p^{+}(s)\|^{2}_{2,\varrho} +\|g(s)\|^{2}_{2,\varrho})\,ds\right).
\end{eqnarray}
\end{stw}
\begin{dow}
Let $\xi_{n}\in C_{0}^{\infty}(\mathbb{R}^{d})$ be a function such
that $\xi_n=1$ on $B(0,n)$. By proposition \ref{stw2.8}, (H2) and
(\ref{eq1.1})
\begin{eqnarray*}
&&\|u(t)\xi_{n}\varrho\|^{2}_{2}+\int_t^T\langle a(s)\nabla
u(s),\nabla (u\xi^{2}_{n}\varrho^{2})(s) \rangle_2\, ds\\
&&\quad=\|\varphi \xi_{n}\varrho\|^{2}_{2}+2\int_t^T \langle
f_u(s),u(s) \xi^{2}_{n}\varrho^{2} \rangle_{2}\,ds
+2\int_t^T\!\!\int_{\mathbb{R}^{d}}u\xi^{2}_{n}\varrho^{2}\,d\mu \\
&&\quad\leq \|\varphi
\xi_{n}\varrho\|^{2}_{2}+\int^T_t(\|g(s)\xi_{n}\varrho\|^{2}_2
+C\|u(s)\xi_{n}\varrho\|^{2}_2
+\frac{\lambda}{2}\|\nabla u(s)\xi_{n}\varrho\|^{2}_2)\,ds\\
&&\qquad+\int_t^T\!\!\int_{\mathbb{R}^{d}}p^{+}\xi^{2}_{n}\varrho^{2}\,d\mu,
\end{eqnarray*}
Moreover, by (\ref{eq2.1}) with $\eta=p^+ \xi^{2}_{n}\varrho^{2}$
we have
\begin{eqnarray*}
\int_t^T\!\!\int_{\mathbb{R}^{d}}p^{+}\xi^{2}_{n}\varrho^{2}\,d\mu
&\leq&
\frac12(\|\varphi\xi_{n}\varrho\|^{2}_{2}+\|p^{+}(T)\xi_{n}\varrho\|^{2}_{2}
+\|u(t)\xi_{n}\varrho\|^{2}_{2}+\|p^{+}(t)\xi_{n}\varrho\|^{2}_{2})\\
&&+\int^T_t(\|g(s)\xi_{n}\varrho\|^{2}_2+\|u(s)\xi_{n}\varrho\|^2_2
+C\|p^{+}(s)\xi_{n}\varrho\|^{2}_2)\,ds\\
&&+\int_t^T(\|\frac{\partial p^{+}}{\partial
s}(s)\xi_{n}\varrho\|^{2}_2) +|\langle a(s)\nabla u(s),\nabla
(p^{+}\xi^{2}_{n}\varrho^{2})(s) \rangle_2|)\,ds
\end{eqnarray*}
By the above estimates and the fact that $|\nabla\varrho|\leq
2\alpha\varrho$ there is $C$ such that
\begin{eqnarray*}
&&\|u(t)\xi_{n}\varrho\|^{2}_{2}+\int^T_t\|\nabla
u(s)\xi_{n}\varrho\|_2^{2}\leq
C\left(\|\varphi\xi_{n}\varrho\|^{2}_{2}
+\sup_{t\in[0,T]}\|p^{+}(t)\xi_{n}\varrho\|^{2}_{2}\right.\\
&&\qquad+ \left.\int^T_t (\|\frac{\partial p^{+}}{\partial
s}(s)\xi_{n}\varrho\|^{2}_2 +\|\nabla p^{+}(s)\xi_{n}\varrho\|^{2}_2
+\|g(s)\xi_{n}\varrho\|^{2}_2)\,ds\right.\\
&&\qquad+\left.\int^T_t\|u(s)\xi_{n}\varrho\|^{2}_2\,ds
+\varepsilon_{n}^{1}+\varepsilon_{n}^{2}\right),
\end{eqnarray*}
where
\begin{eqnarray*}
&&\varepsilon^{1}_{n}=\int^T_t |\langle a(s)\nabla
u(s),u(s)\varrho^2\nabla\xi^2_n\rangle_2|\,ds,\,\,\,
\varepsilon^{2}_{n}=\int^T_t |\langle a(s)\nabla
u(s),p^+(s)\varrho^2\nabla\xi^2_n\rangle_2|\,ds.
\end{eqnarray*}
Since $\varepsilon^{1}_{n}\rightarrow 0$,
$\varepsilon^{2}_{n}\rightarrow 0$ as $n\rightarrow\infty$,
applying Gronwall's lemma we see from the above estimates that
$\sup_{t\in[0,T]}\|u(t)\|^{2}_{2,\varrho}+\int_0^T\|\nabla
u(t)\|^{2}_{2,\varrho}\,dt$ is bounded by the right-hand side of
(\ref{nierownosc1trzyStrony}), which when combined with
(\ref{eq2.06}) and (\ref{eq2.1}) gives
(\ref{nierownosc1trzyStrony}).
\end{dow}
\medskip

For convenience of the reader we now recall  definitions of
solutions of an obstacle problem in the sense of variational
inequalities (see, e.g., \cite{BensoussanLions,Brezis,Lions}).
\medskip\\
{\em Definition} We say that $u$ is a weak solution of
OP$(\varphi,f,h)$ in the variational sense if $u\in
W^{0,1}_{2,\varrho}(Q_T)$ and for any $v\in
W^{1,1}_{2,\varrho}(Q_T)$ such that $v\geq h$,
\begin{eqnarray}\label{dod nierownosc}
&&\int_0^T\langle \frac{\partial v}{\partial
t}(t),(v-u)(t)\rangle_{2,\varrho}\,dt +\int_0^T\langle A_t
u(t),(v-u)(t)\rangle _{2,\varrho}\,dt\nonumber\\
&&\qquad+\int_0^T\langle f_u(t),(v-u)(t)
\rangle_{2,\varrho}\,dt\le
\frac12\|\varphi-v(T)\|^{2}_{2,\varrho}\,,
\end{eqnarray}
where
\[
\langle A_t u(t),(v-u)(t)\rangle _{2,\varrho}=-\frac12\langle
a(t)\nabla u(t),\nabla((v-u)(t)\varrho^{2})\rangle_{2}.
\]
{\em Definition} We say that $u$ is a strong solution of
OP$(\varphi,f,h)$ in the variational sense if
$u\in\mathcal{W}_{\varrho}$, $u(T)=\varphi$ and for any $v\in
W^{0,1}_{2,\varrho}(Q_T)$ such that $v\geq h$,
\begin{eqnarray}\label{dod nierownosc2}
&&\int_0^T\langle \frac{\partial u}{\partial
t}(t),(v-u)(t)\rangle\,dt+\int_0^T\langle A_t
u(t),(v-u)(t)\rangle _{2,\varrho}\,dt\nonumber\\
&&\qquad+\int_0^T\langle f_u(t),(v-u)(t)
\rangle_{2,\varrho}\,dt\le0,
\end{eqnarray}
where $\langle\cdot,\cdot\rangle$ denote the duality pairing
between $W^{1}_{2,\varrho}(\mathbb{R}^{d})$ and
$W^{-1}_{2,\varrho}(\mathbb{R}^{d})$.
\medskip

The following proposition shows that continuous solutions of the
problem (\ref{eq1.4}) coincide with solutions of (\ref{eq1.3}) in
the variational sense.
\begin{stw}
If $(u,\mu)$ is a solution of OP$(\varphi,f,h)$ such that $u\in
W^{0,1}_{2,\varrho}(Q_T)\cap C(\breve Q_T)$ then $u$ is a weak
solution of the problem in the variational sense. If, in addition,
$u\in\mathcal{W}_\varrho$, then  $u$ is a strong solution of
OP$(\varphi,f,h)$ in the variational sense.
\end{stw}
\begin{dow}
Let $u\in W^{0,1}_{2,\varrho}(Q_T)\cap C(\breve Q_T)$ and let
$(u,\mu)$ be a solution of OP$(\varphi,f,h)$. By proposition
\ref{stw2.8},
\begin{eqnarray}
\label{eq2.19} &&\frac12\|u(0)\|^2_{2,\varrho}-\int_0^T\langle
A_tu(t),u(t)\rangle_{2,\varrho}\,dt \nonumber\\
&&\quad=\frac12\|\varphi\|^2_{2,\varrho}+\int_0^T\langle
f_u(t),u(t)\rangle_{2,\varrho}\,dt+\int_{Q_T}u\varrho^2\,d\mu.
\end{eqnarray}
On the other hand, from definition of solution of of OP$(\varphi,f,h)$ it follows that
for any $v\in W^{1,1}_{2,\varrho}(Q_T)$ we have
\begin{eqnarray}
\label{gwiazdka4} &&\int_0^T\langle u(t),\frac{\partial v}{\partial
t}(t)\rangle_{2,\varrho}\,dt -\int_0^T\langle
A_tu(t),v(t)\rangle_{2,\varrho}\,dt\nonumber \\
&&\quad=\int_0^T\langle f_u(t),v(t)\rangle_{2,\varrho}\,dt
+\int_{Q_T}v\varrho^2\,d\mu
+\langle\varphi,v(T)\rangle_{2,\varrho}-\langle
u(0),v(0)\rangle_{2,\varrho}.
\end{eqnarray}
Combining (\ref{eq2.19}) with (\ref{gwiazdka4}) we get
\begin{eqnarray}\label{dod rownosc3}
&&\int_0^T\langle \frac{\partial v}{\partial
t}(t),(v-u)(t)\rangle_{2,\varrho}\,dt +\int_0^T\langle
A_tu(t),(v-u)(t)\rangle_{2,\varrho}\,dt \nonumber\\
&&\qquad+\int_0^T\langle
f_u(t),(v-u)(t)\rangle_{2,\varrho}\,dt\nonumber\\
&&\quad=-\frac12\|u(0)\|_{2,\varrho}^2+\frac12\|\varphi\|_{2,\varrho}^2
+\frac12\|v(T)\|_{2,\varrho}^2-\frac12\|v(0)\|_{2,\varrho}^2
-\langle \varphi,v(T)\rangle_{2,\varrho}\nonumber\\
&&\qquad+\langle u(0),v(0)\rangle_{2,\varrho}
+\int_{Q_T}(u-v)\varrho^2\,d\mu\nonumber\\
&&\quad= -\frac12\|u(0)-v(0)\|_{2,\varrho}^2
+\frac12\|\varphi-v(T)\|_{2,\varrho}^2
+\int_{Q_T}(u-v)\varrho^2\,d\mu.
\end{eqnarray}
Since $v\ge h$, $\int_{Q_T}(u-v)\varrho^2d\,\mu\le 0$, so
(\ref{dod nierownosc}) follows. Now, assume additionally that
$u\in\mathcal{W}_{\varrho}$. Then by (\ref{gwiazdka4}) we have
\begin{eqnarray}
\label{gwiazdka5} &&-\int_0^T\langle\frac{\partial u}{\partial
t}(t),v(t)\rangle\,dt -\int_0^T\langle
A_tu(t),v(t)\rangle_{2,\varrho}\,dt\nonumber \\
&&\quad=\int_0^T\langle f_u(t),v(t)\rangle_{2,\varrho}\,dt
+\int_{Q_T}v\varrho^2\,d\mu
\end{eqnarray}
for every $v\in W^{1,1}_{2,\varrho}(Q_T)$. Let $E\subset\subset
\breve{Q}_{T}$. From (\ref{gwiazdka5}) with a positive $v\in
C^{\infty}_0(\breve{Q}_T)$ such that $v\ge\mathbf{1}_E$ we
conclude that
\[
\mu(E)\le C(\overline{cap}_{\breve{Q}_{T}}(E))^{1/2}
(\|\frac{\partial u}{\partial t}
\|_{\mathbb{L}^{2}([0,T],W^{-1}_{\varrho}(\mathbb{R}^{d}))}
+\|f_{u}\|_{2,\varrho,T}+\|\nabla u\|_{2,\varrho,T}),
\]
where
\[
\overline{cap}_{\breve{Q}_{T}}(E)
=\inf\{\int_{Q_{T}}|\nabla\eta(t,x)|^{2}\,dt\,dx: \eta\in
C_{0}^{\infty}(\breve{Q}_{T}), \eta\geq \mathbf{1}_{E}\}.
\]
On the other hand, it is known (see \cite{B.Mo.}) that
$\overline{cap}_{\,\breve Q_{T}}(E)=\int_{0}^{T}
cap_{\,\mathbb{R}^d}(E_{t})\,dt$. Therefore, if $v\in
W_{2,\varrho}^{0,1}(Q_{T})$, then there is a version of it which
is defined $\overline{q.e.}$. Since we know already that
$\mu\ll\overline{cap}_{\breve{Q}_{T}}$, the integral
$\int_{Q_T}v\,d\mu$ is well defined for $v\in
W_{2,\varrho}^{0,1}(Q_{T})$. Therefore, by approximation argument,
we may take as a test function in (\ref{gwiazdka5}) any $v\in
W_{2,\varrho}^{0,1}(Q_T)$. Now from (\ref{gwiazdka5}) we conclude
that for any $v\in W^{0,1}_{2,\varrho}(Q_T)$ such that $v\geq h$,
\begin{eqnarray*}
&&\int_0^T\langle \frac{\partial u}{\partial
t}(t),(v-u)(t)\rangle_{2,\varrho}\,dt + \int_0^T\langle
A_tu(t),(v-u)(t)\rangle_{2,\varrho}\,dt\\
&&\qquad+\int_0^T\langle f_u(t),(v-u)(t)\rangle_{2,\varrho}\,dt
=\int_{Q_T}(u-v)\varrho^2\,d\mu\le 0,
\end{eqnarray*}
and the proof is complete.
\end{dow}
\medskip

Let us note here that in Theorem \ref{th3.16} we will prove that
if $(u,\mu)$ is a solution of an obstacle problem, then $u$ is the
minimal solution of the same problem in the variational sense.

\nsubsection{Existence, uniqueness and stochastic representation
of solutions of an obstacle problem} \label{sec3}

We begin with a general uniqueness result for continuous solutions
of (\ref{eq1.4}) satisfying some weak integrability assumptions.
\begin{tw}
\label{jednoznacznosc} Assume (H1)--(H3). Then OP$(\varphi,f,h)$ has
at most one solution $(u,\mu)$ such that $u\in C(\breve{Q}_{T})\cap
W^{0,1}_{2,\varrho}(Q_T)$.
\end{tw}
\begin{dow}
Suppose that $(u_{1},\mu_{1}),\ (u_{2},\mu_{2})$ are solutions of
OP$(\varphi,f,h)$ such that $u_{1}, u_{2}\in C(\breve{Q}_{T})\cap
W^{0,1}_{2,\varrho}(Q_T)$ and let $u=u_1-u_2$, $\mu=\mu_1-\mu_2$.
Let $\xi_n:\mathbb{R}^d\rightarrow[0,1]$, $n\in\mathbb{N}$ be a
smooth function such that $\xi_n(x)=1$ if $|x|\le n$ and
$\xi_n(x)=0$ if $|x|\ge n+1$. By the definition of solution of
OP$(\varphi,f,h)$, for any $\eta\in W^{1,1}_{2,0}(Q_T)$ we have
\begin{eqnarray}\label{rownanie214}
&&\langle u(t),\eta(t)\rangle_2+\int^T_t\langle u(s),\frac{\partial
\eta}{\partial s}(s)\rangle_2\,ds+\frac12\int_t^T\langle a(s)\nabla
u(s),\nabla\eta(s)\rangle_2\,ds\nonumber\\
&&\quad=\int^T_t\!\!\int_{\mathbb{R}^d}\eta\,d\mu+\int^T_t\langle
f_{u_1}(s)-f_{u_2}(s),\eta(s)\rangle_2\,ds,\quad t\in[0,T].
\end{eqnarray}
From proposition \ref{stw2.8} we conclude that
\begin{eqnarray}
\label{eq2.9} &&\|u(t)\varrho\xi_n\|^{2}_{2}+\int^T_t\langle
a(s)\nabla
u(s),\nabla(u\varrho^2\xi_n^2)(s)\rangle_2\,ds\nonumber\\
&&\quad=2\int^T_t\!\!\int_{\mathbb{R}^d}u\varrho^2\xi^2_nd\mu
+2\int_{t}^{T}\langle
f_{u_{1}}(s)-f_{u_{2}}(s),u(s)\varrho^{2}\xi^2_n\rangle_{2}\,ds\nonumber\\
&&\quad\le2\int_{t}^{T}\langle
f_{u_{1}}(s)-f_{u_{2}}(s),u(s)\varrho^{2}\xi^2_n\rangle_{2}\,ds,
\end{eqnarray}
the last inequality being a consequence of the fact that
\begin{eqnarray}
\label{eq3.03}
\int^T_t\!\!\int_{\mathbb{R}^d}u\varrho^{2}\xi^2_n\,d\mu
&=&\int^T_t\!\!\int_{\mathbb{R}^d}u_1\varrho^{2}\xi^2_n\,d\mu_{1}
-\int^T_t\!\!\int_{\mathbb{R}^d}u_1\varrho^{2}\xi\,d\mu_{2} \\
&&-\int^T_t\!\!\int_{\mathbb{R}^d}u_2
\varrho^{2}\xi^2_n\,d\mu_{1}+\int^T_t\!\!\int_{\mathbb{R}^d}u_2
\varrho^{2}\xi^2_n\,d\mu_{2}\nonumber \\
&=&\int^T_t\!\!\int_{\mathbb{R}^d}
\varrho^{2}\xi^2_n(u_{1}-h)\,d\mu_{1}
+\int^T_t\!\!\int_{\mathbb{R}^d}
\varrho^{2}\xi^2_n(u_{2}-h)\,d\mu_{2}\nonumber\\
&&+\int^T_t\!\!\int_{\mathbb{R}^d}
\varrho^{2}\xi^2_n(h-u_{2})\,d\mu_{1}
+\int^T_t\!\!\int_{\mathbb{R}^d}
\varrho^{2}\xi^2_n(h-u_{1})\,d\mu_{2}\leq0.\nonumber
\end{eqnarray}
By (\ref{eq2.9}) and (H2),
\begin{eqnarray*}
&&\|u(t)\varrho\xi_n\|^{2}_{2}
+\lambda\int^T_t\|\nabla u(s)\varrho\xi_n\|^2_2\,ds \\
&&\quad\le-\int^T_t \langle a(s)\nabla
u(s),u(s)\xi^2_n\nabla\varrho^2\rangle_2\,ds-\int^T_t \langle
a(s)\nabla
u(s),u(s)\varrho^2\nabla\xi^2_n\rangle_2\,ds \\
&&\qquad+\frac{\lambda}2\int^T_t\|\nabla u(s)\varrho\xi_n\|^2_2\,ds
+2(L+\lambda^{-1}L^2)\int^T_t\|u(s)\varrho\xi_n\|^2_2\,ds.
\end{eqnarray*}
Since $|\nabla\varrho^2|\le2\alpha\varrho^2$, we have
\[
\int^T_t |\langle a(s)\nabla
u(s),u(s)\xi^2_n\nabla\varrho^2\rangle_2|\,ds
\le\frac{\lambda}2\int^T_t\|\nabla u(s)\varrho\xi_n\|^2_2\,ds+
\frac{\Lambda\alpha^{2}}{2\lambda}\int^T_t\|u(s)\varrho\xi_n\|^2_2\,ds.
\]
Consequently, there is $C>0$ not depending on $n$ such that
\[
\|u(t)\varrho\xi_n\|^{2}_{2}\le
C\int^T_t\|u(s)\varrho\xi_n\|^2_2\,ds +\int^T_t|\langle a(s)\nabla
u(s),u(s)\varrho^2_n\nabla\xi_n^2\rangle_2|\,ds
\]
for $t\in[0,T]$. Letting $n\rightarrow\infty$ we get
\[
\|u(t)\|^2_{2,\varrho}\le C\int^T_t\|u(s)\|^2_{2,\varrho}\,ds, \quad
t\in[0,T]
\]
and hence, by Gronwall's lemma, $u=0$, i.e. $u_{1}=u_{2}$. Using
this and (\ref{rownanie214}) we see that
$\int_{Q_{T}}\eta(s,x)d\mu_{1} =\int_{Q_{T}}\eta(s,x)d\mu_{2}$ for
any $\eta\in W^{1,1}_{2,0}(Q_T)$, which shows that $\mu_1=\mu_2$.
\end{dow}
\medskip

To prove existence of a solution of the problem (\ref{eq1.4}) and
its stochastic representation we have to impose additional
integrability assumptions on $g$ and $h$ to ensure existence of a
solution of RBSDE$(\varphi,f,h)$. The assumptions must guarantee
also continuity of $u$ because we are able to prove uniqueness and
a priori estimates only for continuous weak solutions of
OP$(\varphi,f,h)$. Proposition \ref{stw2.3} and Theorem
\ref{tw2.4} therefore suggest that if we want the representation
(\ref{eq1.5}) to hold we should assume at least that
\[
\forall_{K\subset\subset [0,T)\times \mathbb{R}^{d}}\quad
\sup_{(s,x)\in K}(E_{s,x}\sup_{s\leq t\leq T}
|h^{+}(t,X_t)|^2+E_{s,x}\int_s^T|g(t,X_t)|^2\,dt)<\infty.
\]
Our assumptions on $h$ are slightly stronger but nevertheless
seems to be close to the best possible.

Now we provide a useful inequality for moments of the diffusion
$(X,P_{s,x})$. It is perhaps known but we could not find a proper
reference. The inequality is given only for moments greater or
equal to 4, because such a form is sufficient for our purposes.

\begin{stw}\label{doob}
If $(X,P_{s,x})$ is a Markov process associated with $A_t$ then
for every $p\geq 4$,
\[
E_{s,x}\sup_{s\leq t \leq T}|X_{t}|^{p}\leq CE_{s,x}|X_{T}|^{p},
\]
where $C$ depends only on $\lambda,\Lambda,d$ and $T$.
\end{stw}
\begin{dow}
Let $u_{n}$ be a solution of PDE$(\varphi_n,0)$ with $\varphi_{n}(x)=
|x|^{p/2}\mathbf{1}_{B(0,n)}(x)$. From \cite{A.Roz.BSDE} we know
that the pair $(u_{n}(t,X_{t}),\sigma\nabla u_{n}(t,X_{t}))$,
$t\in[s,T]$, is a solution of BSDE$(\varphi_n,0)$, i.e
\[
u_{n}(t,X_{t})=\varphi_{n}(X_{T}) -\int_{t}^{T}\sigma\nabla
u_{n}(\theta,X_{\theta})\,dB_{s,\theta},\quad P_{s,x}\mbox{-}a.s.
\]
from which we obtain in particular that
$u_{n}(s,x)=E_{s,x}\varphi_{n}(X_{T})$. It is known that
$u_n\rightarrow u$ uniformly in compact subsets of $Q_T$. By
Aronson's lower estimate, for all sufficiently large
$n\in\mathbb{N}$ we have
\begin{eqnarray*}
|u_{n}(s,x)|&=&E_{s,x}|X_{T}|^{p/2}\mathbf{1}_{B(0,n)}(X_{T})\\
&\ge& C\int_{B(0,n)}|y|^{p/2}(T-s)^{-d/2}
\exp(-\frac{|y-x|^{2}}{C(T-s)})\,dy\\
&=&C^{1+(d/2)}
E\mathbf{1}_{B(0,n)}(X_{C(T-s)}+x)|X_{C(T-s)}+x|^{p/2}\\
&\ge& C^{1+(d/2)}\left(E\mathbf{1}_{B(0,n)}(X_{C(T-s)}+x)
(|x|^{2}+2\sum^d_{i=1}x_iX^i_{C(T-s)})\right)^{p/4}
\end{eqnarray*}
($E$ denotes expectation with respect to the standard Wiener
measure on $\Omega$). Letting $n\rightarrow\infty$ we see that
$|u(s,x)|\ge C^{1+(d/2)}|x|^{p/2}$. By the above and known a
priori estimates for BSDE we get
\begin{eqnarray*}
E_{s,x}\sup_{s\le t\le T}|X_t|^p&\le& CE_{s,x}\sup_{s\le t\le T}
|u(t,X_t)|^2\le C\liminf_{n\rightarrow\infty}E_{s,x}\sup_{s\le
t\le T} |u_n(t,X_t)|^2\\
&\le& C\liminf_{n\rightarrow\infty}E_{s,x} |\varphi_n(X_T)|^2\le
CE_{s,x}|X_T|^{p},
\end{eqnarray*}
which completes the proof.
\end{dow}
\medskip

Here and subsequently, we write $\mu_n\Rightsx\mu$ if for fixed
$(s,x)\in[0,T)\times\mathbb{R}^d$,
\[
\int_{Q_T}\xi(t,y)p(s,x,t,y)\,d\mu_n(t,y)\rightarrow
\int_{Q_T}\xi(t,y)p(s,x,t,y)\,d\mu(t,y)
\]
for every $\xi\in C_b(Q_T)$. We use the symbol
``$\Rightarrow$" to denote weak convergence of measures.

\begin{lm}\label{lm3.4}
Let $S$ be a Polish space and let $\mu$, $\mu_{n}$,
$n\in\mathbb{N}$, be probability measures on $S$ such that
$\mu_n\Rightarrow\mu$. If $f,f_n:S\rightarrow\mathbb{R}$ are
continuous functions such that $f_{n}\rightarrow f$ uniformly in
compact subsets of $S$ and
\[
\sup_{n\ge1}\int_{S}|f_{n}|\,d\mu_{n}<\infty,\quad
\lim_{\alpha\rightarrow \infty}
\sup_{n\ge1}\int_{S}|f_{n}|\mathbf{1}_{\{|f_{n}|\ge
\alpha\}}\,d\mu_{n}=0
\]
then
\[
\int_{S}f_{n}\,d\mu_{n}\rightarrow \int_{S}f\,d\mu.
\]
\end{lm}
\begin{dow}
It is sufficient to modify slightly the proof of \cite[Lemma
8.4.3]{Bo}. We omit the details.
\end{dow}
\medskip

We now prove our main existence and representation results. For
reasons to be explained later on, we decided to consider
separately the case of square-integrable data $\varphi,g,h$ and
the case where the data are square-integrable with some weight
$\varrho\in W$ such that $\varrho<1$.

\begin{tw}
\label{tw3.4} Let assumptions (H1)--(H3) hold with
$\varphi\in\mathbb{L}_2(\mathbb{R}^d)$, $g\in\mathbb{L}_2(Q_T)$
and moreover, assume that $g$ satisfies (\ref{n2.2}),
$h\in\mathbb{L}_2(Q_{T})\cap C(Q_{T})$, $h\le\psi$ for some $\psi$
such that $\psi\in W^{1,1}_2(Q_{T})$ and $h\leq c\bar\varrho^{-1}$
for some $c>0$, $\bar\varrho\in W$. Then there exists a unique
weak solution $(u,\mu)$ of OP$(\varphi,f,h)$ such that
\begin{enumerate}
\item[\rm(i)]
$u\in C([0,T)\times\mathbb{R}^{d})$,
\item[\rm(ii)]
$u_{n}\uparrow u$, $u_{n}\rightarrow u$ in
$W^{0,1}_{2,loc}(Q_{T})\cap C([0,T],\mathbb{L}_{2}^{loc}(Q_{T}))$,
$u_n\rightarrow u$ in $\mathbb{L}_2(Q_T)$ and
$\mu_{n}\Rightarrow\mu$, $\mu_{n}\rightarrow\mu$ in
$(W_{2,loc}^{1,1}(Q_{T}))^{*}$, $\mu_n\Rightsx\mu$ for every
$(s,x)\in[0,T)\times\mathbb{R}^d$, where
$d\mu_{n}=n(u_{n}-h)^{-}d\lambda$ and $u_n$ is a unique weak
solution of the Cauchy problem
\begin{equation}
\label{eq3.2} (\frac{\partial}{\partial t}+A_t)
u_{n}=-f_{u_{n}}-n(u_{n}-h)^{-},\quad u_{n}(T)=\varphi.
\end{equation}
\end{enumerate}
Moreover, for each $(s,x)\in[0,T)\times{\mathbb{R}^d}$,
\begin{equation}
\label{eq3.02} (u(t,X_{t}),\sigma \nabla
u(t,X_{t}))=(Y^{s,x}_{t},Z^{s,x}_{t}),\quad t\in[s,T],\quad
P_{s,x}\mbox{-}a.s.
\end{equation}
and
\begin{equation}
\label{eq3.04} E_{s,x}\int_{s}^{T}\xi(t,X_t)\,dK^{s,x}_t
=\int_s^T\!\!\int_{\mathbb{R}^d}\xi(t,y)p(s,x,t,y)\,d\mu(t,y)
\end{equation}
for every $\xi\in C_{b}(Q_{T})$, where $(Y^{s,x},Z^{s,x},K^{s,x})$
is a solution of RBSDE$(\varphi,f,h)$.
\end{tw}
\begin{dow} {\em Step 1.} We first show existence of $u\in W^{0,1}_2(Q_T)$
and a Radon measure $\mu$ on $Q_{T}$ such that
\begin{eqnarray}
\label{eq3.5} &&\int_{0}^{T}\langle
u(t),\frac{\partial\eta}{\partial t}(t)\rangle_{2}\,dt
+\frac12\int_{0}^{T}\langle a(t)\nabla
u(t),\nabla\eta(t)\rangle_{2}\,dt \nonumber\\
&&\qquad=\int_{0}^{T}\langle f_{u}(t),\eta(t)\rangle_{2}\,dt
+\int_{Q_T}\eta\,d\mu +\langle\varphi,\eta(T)\rangle_{2}
\end{eqnarray}
for every $\eta\in W^{1,1}_{2}(Q_{T})\cap C(Q_{T})$ such that
$\eta(0)\equiv0$.  From Proposition \ref{stw2.3} we know that
there exists a unique weak solution $u_n$ of (\ref{eq3.2}) such
that $u_{n}\in C([0,T],\mathbb{L}_{2}(\mathbb{R}^{d}))$ and
$u_n\in C([0,T)\times\mathbb{R}^d)$. Set $r_{n}=n(u_{n}-h)^{-}$
and let $d\mu_{n}=r_nd\lambda$. Then for any $\eta\in
W^{1,1}_{2,0}(Q_{T})$,
\begin{eqnarray}
\label{eq3.4} &&\int_{t}^{T}\langle
u_{n}(s),\frac{\partial\eta}{\partial s}(s)\rangle_{2}\,ds
+\frac12\int_{t}^{T}\langle a(s)\nabla u_{n}(s),\nabla
\eta(s)\rangle_{2}\,ds =\int_{t}^{T}\langle
f_{u_{n}}(s),\eta(s)\rangle_{2}\,ds \nonumber\\
&&\qquad +\int_t^T\!\!\int_{\mathbb{R}^d}\eta\,d\mu_{n}
+\langle\varphi,\eta(T)\rangle_{2} -\langle
u_{n}(t)\eta(t)\rangle_{2},\quad t\in[0,T].
\end{eqnarray}
By Proposition \ref{stw2.30} with $p=\psi$ there is $C>0$ such
that
\begin{equation}
\label{eq3.3} \sup_{t\in[0,T]}\|u_{n}(t)\|^{2}_{2}
+\int^T_0|\langle a(t)\nabla u_{n}(t),\nabla
u_{n}(t)\rangle_2|\,dt \leq C
\end{equation}
for every $n\in\mathbb{N}$. Since, by continuity of $u_{n}$ and
comparison results (see \cite[Theorem 4.1.7]{BensoussanLions}),
$u_{n}(t,x)\leq u_{n+1}(t,x)$ for every $(t,x)\in
[0,T)\times\mathbb{R}^d$, there is $u$ such that $u_{n}\uparrow u$.
By Fatou's lemma and (\ref{eq3.3}), $u\in\mathbb{L}_{2}(Q_T)$. In
fact, since $u_1\le u_n\le u$, it follows that $u^2_n\le u_1^2+u^2$
and hence, by the Lebesgue dominated convergence theorem, that
$u_n\rightarrow u$ in $\mathbb{L}_{2}(Q_T)$. Let
$f_n(t,x,y)=f_{u_n}(t,x)+n(y-h(t,x))$ and let
$(Y^{s,x,n},Z^{s,x,n})$ be a unique solution of BSDE$(\varphi,f_n)$.
By results of \cite{A.Roz.BSDE}, $(u_n(t,X_{t}),\sigma\nabla
u_n(t,X_{t}))=(Y^{s,x,n}_t,Z^{s,x,n}_{t})$ $P_{s,x}$-a.s. and hence,
by (\ref{eq2.4}),
\begin{eqnarray*}
|u_{n}(s,x)|^{2}&\leq& CE_{s,x}\left(|\varphi(X_T)|^2
+\int_{s}^{T}|g(t,X_t)|^2\,dt +\sup_{s\leq t\leq T}
|h^{+}(t,X_{t})|^{2}\right)
\end{eqnarray*}
for all $(s,x)\in[0,T)\times\mathbb{R}^d$. Thus,
$u(s,x)=\sup_{n\ge1}u_{n}(s,x)<\infty$ and consequently, $u$ is
lower semi-continuous on $[0,T)\times\mathbb{R}^d$. Let
$K\subset\subset Q_{T}$ and let $\eta\in
C^{\infty}_{0}(\breve{Q}_{T})$ be a positive function such that
$\eta=1$ on $K$. Since
\begin{eqnarray*}
&&\mu_{n}(K)\leq \int_{Q_{T}}\eta\,d\mu_{n}=\int_{0}^{T}\langle
u_n(s),\frac{\partial\eta}{\partial
s}(s)\rangle_{2}\,ds+\int_{0}^{T}\langle a(s)\nabla
u_n(s),\nabla\eta(s)\rangle_{2}\,ds\\
&&\qquad\qquad\qquad\qquad\qquad-\int_{0}^{T}\langle
f_{u_{n}}(s),\eta(s)\rangle_{2}\,ds,
\end{eqnarray*}
we conclude from (\ref{eq3.3}) that
$\sup_{n\ge1}\mu_{n}(K)<\infty$. Thus, by the weak compactness
theorem for measures  (see Section 1.9 in \cite{Evance}),
$\{\mu_n\}$ is tight. Therefore there is a subsequence, still
denoted by $\{n\}$, such that $\int_{Q_{T}} f\,d\mu_{n}\rightarrow
\int_{Q_{T}} f\,d\mu$ for every $f\in C_0(Q_{T})$. Let $\varrho\in
W$ be such that
$\int_{\mathbb{R}^d}(\varrho\bar\varrho^{-1}(x))^2\,dx<\infty$ and
let $K\subset\subset Q_T$. Then by Theorem \ref{tw2.2},
\begin{equation}
\label{eq3.010} \int_{\mathbb{R}^{d}}\left(E_{0,x}\int_{0}^{T}
|\nabla(u_{n}-u_m)(t,X_{t})|^2\,dt\right)\varrho^2(x)\,dx \ge
C\|\nabla(u_{n}-u_{m})\varrho\|^{2}_{2}.
\end{equation}
By (\ref{zbieznosci}), for every $x\in\mathbb{R}^d$,
$\xi_{n,m}(x)\equiv E_{0,x}\int_{0}^{T}
|\nabla(u_{n}-u_m)(t,X_{t})|^2\,dt\rightarrow0$ as
$n,m\rightarrow\infty$. Moreover, by (\ref{eq2.4}) and Proposition
\ref{doob},
\[
|\xi_{n,m}(x)|\le C(E_{0,x}|\varphi(X_T)|^2
+E_{0,x}\int^T_0|g(t,X_t)|^2\,dt+\bar\varrho^{-1}(x))
\]
for some $C$ not depending on $n,m$. Therefore it follows from
Theorem \ref{tw2.2} and the Lebesgue dominated convergence theorem
that the left-hand side of (\ref{eq3.010}) converges to zero as
$n,m\rightarrow\infty$ and hence that
$\|\mathbf{1}_K\nabla(u_{n}-u_{m})\|^{2}_2\rightarrow 0$ for any
$K\subset\subset Q_{T}$. Using properties of $\{u_n\}$ and
$\{\mu_n\}$ we have already proved we conclude from (\ref{eq3.4})
that (\ref{eq3.5}) holds for every $\eta\in W^{1,1}_{2}(Q_{T})\cap
C_{0}(Q_{T})$ such that $\eta(0)=0$.
\medskip\\
{\em Step 2.} $u\in C([0,T)\times{\mathbb{R}}^d)\cap
C([0,T],\mathbb{L}^{loc}_2({\mathbb{R}}^d))$. To see this, we first
observe that $u(s,x)=Y^{s,x}_s$ for
$(s,x)\in[0,T)\times{\mathbb{R}^d}$, since
\[
|u(s,x)-Y^{s,x}_{s}|^{2}\le2\lim_{n\rightarrow\infty}
(|(u-u_{n})(s,x)|^2+E_{s,x}|Y^{s,x,n}_{s}-Y^{s,x}_{s}|^{2})=0.
\]
Hence $u(s,x)=Y^{s,x}_s\le h(s,X_s)=h(s,x)$, i.e. $(u-h)^-=0$ and,
by (\ref{ciaglosc}), for any $n,m\in\mathbb{N}$, $\delta>0$ and
any $K\subset\subset[0,T-2\delta]\times{\mathbb{R}}^d$ we have
\begin{equation}
\label{eq3.6} |u_{n}(s,x)-u_{m}(s,x)|^{2} \le
C(E_{s,x}|Y^{s,x,n}_{T-\delta}-Y^{s,x,m}_{T-\delta}|^{2}
+I^{s,x}_{n,m}+I^{s,x}_{m,n}),
\end{equation}
where
\[
I^{s,x}_{n,m}=E_{s,x}\int_{s}^{T-\delta}
(Y^{s,x,n}_{t}-h(t,X_{t}))^{-}\,dK_{t}^{s,x,m}
\]
and $K^{s,x,m}$ is defined as in Theorem \ref{tw2.4}. By Aronson's
upper estimate,
\begin{eqnarray*}
E_{s,x}|Y^{s,x,n}_{T-\delta}-Y^{s,x,m}_{T-\delta}|^{2}
&=&E_{s,x}|(u_{n}-u_{m})(T-\delta,X_{T-\delta})|^{2}\\
&=&\int_{\mathbb{R}^{d}}|(u_{n}-u_{m})(T-\delta,y)|^{2}
p(s,x,T-\delta,y)\,dy \\
&\leq& C\delta^{-d/2}\|u_{n}-u_{m}\|^{2}_{2}
\end{eqnarray*}
with some $C$ depending neither on $(s,x)\in K$ nor on
$n,m\in\mathbb{N}$. Moreover,
\begin{eqnarray*}
|I^{s,x}_{n,m}|^2\le E_{s,x}\sup_{s\leq t\leq
T-\delta}|(u_{n}(t,X_{t})-h(t,X_{t}))^{-}|^{2} \cdot
E|K^{s,x,m}_{T-\delta}|^{2}.
\end{eqnarray*}
In view of (\ref{eq2.4}), $\sup_{n\ge1}\sup_{(s,x)\in
K}E|K^{s,x,n}_{T-\delta}|^{2}<\infty$. By Dini's theorem,
$(u_{n}-h)^{-}\rightarrow(u-h)^{-}=0$ uniformly in any compact
subset of $[0,T)\times{\mathbb{R}}^d$. Therefore, since
$|h(t,X_t)|\le|\bar u(t,X_t)|+\bar\varrho^{-1}(X_t)$, $t\in[0,T]$,
where $\bar u$ is a solution of PDE$(\varphi,f)$, it follows from
Theorem \ref{tw2.1}, Proposition \ref{doob}, Lemma \ref{lm3.4} and
the Lebesgue dominated convergence theorem that $\sup_{(s,x)\in K}
(I^{s,x}_{n,m}+I^{s,x}_{m,n})\rightarrow0$ as
$n,m\rightarrow\infty$. From the above estimates and (\ref{eq3.6})
it follows that $u_{n}\rightarrow u$ uniformly in any compact
subset of $[0,T)\times\mathbb{R}^{d}$, which implies continuity of
$u$ on $[0,T)\times\mathbb{R}^{d}$. By Theorem \ref{tw2.2},
\begin{eqnarray*}
&&\int_{\mathbb{R}^{d}}E_{0,x}\sup_{t\in[0,T]}
|Y^{0,x,n}_{t}-Y^{0,x,m}_{t}|^{2}\varrho^2(x)\,dx\\
&&\qquad\geq\sup_{t\in[0,T]}\int_{\mathbb{R}^{d}}
E_{0,x}|Y^{0,x,n}_{t}-Y^{0,x,m}_{t}|^{2}\varrho^2(x)\,dx\\
&&\qquad=\sup_{t\in[0,T]}\int_{\mathbb{R}^{d}}
E_{0,x}|(u_{n}-u_m)(t,X_{t})|^{2}\varrho^2(x)\,dx
\ge c\sup_{t\in[0,T]}\|(u_n-u_m)(t)\|^{2}_{2,\varrho}\,.
\end{eqnarray*}
Therefore, choosing $\varrho\in W$  such that
$\int_{\mathbb{R}^d}(\varrho\bar\varrho^{-1}(x))^2\,dx<\infty$ and
arguing as in the proof of convergence of the left-hand side of
(\ref{eq3.010}) we deduce from the above that
$\sup_{t\in[0,T]}\|(u_{n}-u_{m})(t)\|^{2}_{2,\varrho}\rightarrow0$
as $n,m\rightarrow\infty$, and hence that $u\in
C([0,T],\mathbb{L}^{loc}_2({\mathbb{R}}^d))$.
\medskip\\
{\em Step 3.} $u$ is the unique weak solution of the problem
OP$(\varphi,f,h)$. We know that
\begin{eqnarray}
\label{eq3.9} E_{0,x}\int_0^T\xi(t,X_t)\,dK_t^{0,x,n}
&=&E_{0,x}\int_0^T\xi n(u_n-h)^{-}(t,X_t)\,dt\nonumber\\
&=&\int_{Q_T}\xi(t,y)p(s,x,t,y)\,d\mu_n(t,y)
\end{eqnarray}
for all $\xi\in C_0(Q_T)$. Let $K\subset\subset \mathbb{R}^{d}$ and  $\{\xi_n\}\subset
C^+_0(Q_T)$ be such that $\xi_n\downarrow\mathbf{1}_{\{0\}\times
K}$. Since $\mu_n\Rightarrow\mu$, it follows from (\ref{eq3.9}) and (\ref{zbieznosci})
that
\[
\int_{\mathbb{R}^{d}}(E_{0,x}\int_0^T\xi_n(t,X_t)\,dK^{0,x}_t)\,dx\ge
\int_{Q_T}\xi_n(t,y)\,d\mu(t,y)
\]
for $n\in\mathbb{N}$. Letting $n\rightarrow\infty$ and taking
into account that  $K^{0,x}$ is continuous we deduce from the
above inequality that $\int_{Q_T}\mathbf{1}_{\{0\}\times
K}(\theta,y)\,d\mu(\theta,y)=0$. Therefore
$\mu(\{0\}\times K)=0$ for any $K\subset\subset \mathbb{R}^{d}$
and hence $\mu(\{0\}\times\mathbb{R}^d)=0$. Now from Lemma
\ref{lemat o rownowaznosci} and {\em Step 2}\, it follows that
$\mu(\{t\}\times\mathbb{R}^d)=0$ for all $t\in[0,T]$.  Using this
and Lemma \ref{lemat o rownowaznosci} we see that (\ref{eq2.1}) is
satisfied for any $\eta\in W_{2,0}^{1,1}(Q_{T})$ and  $t\in
[0,T]$. Since $u_{n}-h\rightarrow u-h$ uniformly in compact
subsets of $[0,T)\times\mathbb{R}^{d}$,
\begin{eqnarray*}
\int_{s}^{T-\delta}\!\!\int_{\mathbb{R}^{d}}
\xi(u_{n}-h)(t,x))\,d\mu_{n}(t,x) \rightarrow\int_{s}^{T-\delta}\!\!
\int_{\mathbb{R}^{d}}\xi(u-h)(t,x))\,d\mu(t,x)
\end{eqnarray*}
for any $\xi\in C^+_0(\mathbb{R}^{d})$. Hence, since
\begin{eqnarray*}
&&\int_{s}^{T-\delta}\!\!\int_{\mathbb{R}^{d}}
\xi(u_{n}-h)(t,x)\,d\mu_{n}(t,x)\\
&&\qquad =n\int_{s}^{T-\delta}\!\!\int_{\mathbb{R}^{d}}
\xi(u_{n}-h)\cdot(u_{n}-h)^{-}(t,x)\,dt\,dx\leq 0
\end{eqnarray*}
and $u\geq h$, it follows that $\int_{Q_T}\xi(u-h)\,d\mu=0$, which
shows that $u$ solves OP$(\varphi,f,h)$. Uniqueness follows from
Theorem \ref{jednoznacznosc}.
\medskip\\
{\em Step 4.} We show (\ref{eq3.02}). From \cite{A.Roz.BSDE} we know
that
\begin{eqnarray}
\label{rownanie penalizacyjne 1} u_{n}(t,X_{t})&=&\varphi(X_{T})
+\int_{t}^{T}f(\theta,X_{\theta},u_{n}(\theta,X_{\theta}),\sigma\nabla
u_{n}(\theta,X_{\theta}))\,d\theta\nonumber\\
&&+K^{s,x,n}_{T}-K^{s,x,n}_{t}-\int_{t}^{T}\langle \sigma\nabla
u_{n}(\theta,X_{\theta}),dB_{s,\theta}\rangle,\quad
P_{s,x}\mbox{-}a.s.
\end{eqnarray}
for all $n\in\mathbb{N}$. Since $u_{n}\rightarrow u$ uniformly on
compact sets in $[0,T)\times\mathbb{R}^d$, it follows from
(\ref{zbieznosci}) that $u(t,X_{t})=Y^{s,x}_{t}$, $t\in[s,T]$,
$P_{s,x}$\,-a.s. To prove that $\sigma\nabla u(t,X_t)=Z^{s,x}_t$,
$\lambda\otimes P_{s,x}$-a.s. observe that for any
$K\subset\subset\mathbb{R}^d$ and any $\delta\in(0,T-s]$,
\begin{eqnarray*}
&&E_{s,x}\int_{s+\delta}^{T}\mathbf{1}_{\{X_t\in K\}}
|\sigma\nabla u(t,X_t)-Z^{s,x}_t|^{2}\,dt \\
&&\quad\le 2E_{s,x}\int_{s+\delta}^{T}\mathbf{1}_{\{X_t\in
K\}}(|\sigma\nabla(u-u_{n})(t,X_t)|^{2}
+|Z^{s,x,n}_t-Z^{s,x}_t|^{2})\,dt\\
&&\quad\le C\delta^{-d/2}\int_{s+\delta}^{T}\!\int_{K}
|\nabla(u-u_{n})(t,y)|^{2}\,dt\,dy
+2E_{s,x}\int_{s+\delta}^{T}|Z^{s,x,n}_t-Z^{s,x}_t|^{2}\,dt,
\end{eqnarray*}
which converges to zero as $n\rightarrow\infty$ since $\nabla
u_n\rightarrow\nabla u$ in $\mathbb{L}^{loc}_{2}(Q_{T})$. Hence, by
Fatou's lemma, $E_{s,x}\int_{s}^{T}|\sigma\nabla
u(t,X_t)-Z^{s,x}_t|^{2}\,dt=0$, as required.
\medskip\\
{\em Step 5.} We show (\ref{eq3.04}) and that $\mu_n\Rightsx\mu$,
$\mu_{n}\rightarrow\mu$ in $(W_{2,loc}^{1,1}(Q_{T}))^{*}$ for
every $(s,x)\in[0,T)\times\mathbb{R}^d$. Let $\xi=\xi^+-\xi^{-}\in
C_{b}(Q_{T})$. By Theorem \ref{tw2.4},
\begin{equation}
\label{eq3.14}
E_{s,x}\left|\int_{t}^{T}n\xi(u_{n}-h)^-(\theta,X_{\theta})
\,d\theta-\int_{t}^{T}\xi(\theta,
X_{\theta})dK^{s,x}_{\theta}\right|^2\rightarrow0
\end{equation}
for every $t\in[s,T]$. Since $\mu_{n}\Rightarrow\mu$ on $Q_T$ and
$\mu(\{t\}\times\mathbb{R}^d)=0$ for every $t\in[0,T]$, we see
that $\mu_n|_{[t_1,t_2]\times\mathbb{R}^d}
\Rightarrow\mu|_{[t_1,t_2]\times\mathbb{R}^d}$ for every $t\le
t_1\le t_2\le T$. Hence we have
\begin{eqnarray*}
E_{s,x}\int_{s+\delta}^{T}\xi^{+}(t,X_{t})\,dK^{s,x}_{\theta}
&=&\lim_{n\rightarrow\infty} E_{s,x}\int_{s+\delta}^{T}
\xi^{+}(t,X_{t})\,dK^{s,x,n}_{t}\\
&=&\lim_{n\rightarrow\infty}
\int_{s+\delta}^{T}\int_{\mathbb{R}^{d}}\xi^{+}
(t,y)p(s,x,t,y)\,d\mu_{n}(t,y)\\
&=&\int_{s+\delta}^{T}\int_{\mathbb{R}^{d}}
\xi^{+}(t,y)p(s,x,t,y)\,d\mu(t,y)
\end{eqnarray*}
for  $0\le s<T$, $\delta>0$. Applying the monotone convergence
theorem we see that
\begin{equation}
\label{eq3.10}
E_{s,x}\int_s^T\xi^{+}(\theta,X_{\theta})\,dK^{s,x}_{\theta}
=\int_{s}^{T}\!\!\int_{\mathbb{R}^{d}}
\xi^{+}(\theta,y)p(s,x,\theta,y)\,d\mu(\theta,y)
\end{equation}
for all $s\in[0,T]$. In the same manner we can see that
(\ref{eq3.10}) holds for $\xi^{-}$ in place of $\xi^+$ and hence
for $\xi$ in place of $\xi^+$. Consequently, (\ref{eq3.04}) holds
for $s\in[0,T]$. That $\mu_n\Rightsx\mu$ now follows from
(\ref{zbieznosci}), (\ref{eq3.14}). Strong  convergence of
$\{\mu_{n}\}$ to $\mu$ in $(W_{2,loc}^{1,1}(Q_{T}))^{*}$ follows
from (\ref{eq3.4}), (\ref{eq2.1}) and the fact that
$u_n\rightarrow u$ in $W^{0,1}_{2,loc}(Q_T)$.
\end{dow}

\begin{wn}
\label{wn3.5} Under the assumption of Theorem \ref{tw3.4}, for any
$0\le t_1<t_2\le T$ and any closed subset $F$ of $\mathbb{R}^d$ we
have
\begin{equation}
\label{eq3.12} \mu([t_1,t_2]\times F)
=\int_{\mathbb{R}^d}E_{t_1,x}\int_{t_1}^{t_2}
\mathbf{1}_{F}(X_t)\,dK_t^{t_1,x}\,dx.
\end{equation}
\end{wn}
\begin{dow}
Let us choose a sequence $\{\xi_n\}\subset C_b(Q_T)$ of positive
functions such that $\xi_n\downarrow \mathbf{1}_{[t_1,t_2]\times
F}$. Since (\ref{eq3.14}) holds for $\xi_n$ in place of $\xi^+$,
we get (\ref{eq3.12}) letting $n\rightarrow\infty$ and then
integrating with respect to the space variable.
\end{dow}

\begin{wn}
\label{wn3.6} Let assumptions of Theorem \ref{tw3.4} hold. Then
$\mu$ is absolutely continuous with respect to the Lebesgue measure
with density $r$ iff
\begin{equation}
\label{eq3.016}
K^{s,x}_t=\int_s^tr(\theta,X_\theta)\,d\theta,\quad
t\in[s,T],\quad P_{s,x}\mbox{-}a.s.
\end{equation}
\end{wn}
\begin{dow}
Sufficiency follows immediately  from (\ref{eq3.04}). To prove
necessity, suppose that $(u,rd\lambda)$ is a weak solution of the
OP$(\varphi,f,h)$ i.e. $(u-h)d\mu=0$, $u\geq h$ and
\begin{equation}
\label{eq3.16} \frac{\partial u}{\partial t}+A_tu=-f_u-r,\quad
u(T)=\varphi.
\end{equation}
Set
$r^\varepsilon=(r\wedge\varepsilon^{-1})\mathbf{1}_{B(0,\varepsilon^{-1})}$
and let $u_{\varepsilon}$ be a weak solution of
PDE$(\varphi,f_u+r^\varepsilon)$, i.e.
\[
\frac{\partial u_{\varepsilon}}{\partial
t}+A_tu_\varepsilon=-f_u-r^{\varepsilon},\quad
u_\varepsilon(T)=\varphi.
\]
Then
\begin{eqnarray*}
\|u_\varepsilon(t)\|_2^2+\|\sigma\nabla u_\varepsilon\|^2_{2,T}
&=&2\int_t^T\langle f_u(s),u_\varepsilon(s)\rangle_{2}\,ds
+2\int^{T}_{t}\langle
r^\varepsilon(s),u_\varepsilon(s)\rangle_{2}\,ds\\
&\le&\int_t^T\|u_\varepsilon(s)\|^2_2\,ds+\|f_u\|^2_{2,T}
+2\int^{T}_{t}\langle r(s),|u(s)|\rangle_{2}\,ds
\end{eqnarray*}
and hence, by Gronwall's lemma,
\begin{eqnarray}
\label{eq3.28} \|u_\varepsilon(t)\|_2^2+\|\sigma\nabla
u_\varepsilon\|^2_{2,T}\le C(\|f_u\|^2_{2,T}+\int^{T}_{0}\langle
r(t),|u(t)|\rangle_{2}\,dt),
\end{eqnarray}
which is bounded, because $\int^{T}_{0}\langle
r(t),|u(t)|\rangle_{2}\,dt<\infty$ by Proposition \ref{stw2.8}.
Since $\{u_\varepsilon\}$ is increasing, there is $\bar u$ such
that $u_\varepsilon\uparrow\bar{u}$. Since we know that
$\{u_{\varepsilon}\}$ is bounded in $W^{0,1}_2(Q_T)$,
$u_\varepsilon\rightarrow\bar{u}$ in $\mathbb{L}_2(Q_T)$ and
$\nabla u_\varepsilon\rightarrow\nabla\bar{u}$ weakly in
$\mathbb{L}_2(Q_T)$ from which it may be concluded that $\bar u$
is a weak solution of (\ref{eq3.16}). Therefore, $u=\bar{u}$, by
uniqueness of solution of PDE$(\varphi,f_{u}+r)$. Now, define
$r_n,\mu_n$ as in Theorem \ref{tw3.4}. Let $\xi\in C_b(Q_T)$.
Since $\mu_n\Rightsx\mu$ for every
$(s,x)\in[0,T)\times\mathbb{R}^d$ and $d\mu=r\,d\lambda$,
\begin{eqnarray}
\label{eq3.19} E_{s,x}\int_s^T\xi(t,X_t)\,dK_t^{s,x}
=E_{s,x}\int_s^T\xi(t,X_t)r(t,X_t)\,dt.
\end{eqnarray}
Indeed, for every $n\in\mathbb{N}$ we have
\begin{eqnarray*}
E_{s,x}\int_s^T\xi(t,X_t)\,dK^{s,x,n}_t
&=&\int_{\mathbb{R}^d}\!\int_s^T\xi(t,y)p(s,x,t,y)\,d\mu_n(t,y)\\
&=&\int_{\mathbb{R}^d}\!\int_s^T\xi(t,y)p(s,x,t,y)r_n(t,y)\,dt\,dy,
\end{eqnarray*}
so letting $n\rightarrow\infty$ leads to (\ref{eq3.19}). By
approximation argument, (\ref{eq3.19}) holds for any  $\xi\in
C(Q_T)$ such that
$E_{s,x}\int_s^T|\xi(t,X_t)|\,dK^{s,x}_t<\infty$. In particular,
it holds for $\xi=u-h$. Hence, letting $\varepsilon\downarrow0$
and applying the Lebesgue dominated convergence theorem we obtain
\begin{eqnarray}
\label{eq3.31} &&I^1_{\varepsilon}\equiv
E_{s,x}\int_s^T(u_\varepsilon-h)(t,X_t)
r^\varepsilon(t,X_t)\,dt\\
&&\qquad\rightarrow E_{s,x}\int_s^T(u-h)(t,X_t)
r(t,X_t)\,dt=E_{s,x}\int_s^T(u-h)(t,X_t)\,dK^{s,x}_t=0.\nonumber
\end{eqnarray}
By representation results proved in  \cite{A.Roz.BSDE},
\begin{eqnarray}\label{gwiazdka}
u_\varepsilon(t,X_t)=\varphi(X_T)
+\int_t^T(f_u+r^{\varepsilon})(\theta,X_\theta)\,d\theta
-\int_t^T\sigma\nabla
u_\varepsilon(\theta,X_\theta)\,dB_{s,\theta}.
\end{eqnarray}
Applying It\^{o}'s formula we obtain
\begin{eqnarray*}
I^2_{\varepsilon}&\equiv&E_{s,x}|u_\varepsilon(t,X_t)|^2
+E_{s,x}\int_t^T|\sigma\nabla
u_\varepsilon(\theta,X_{\theta}^{s,x})|^2\,dt\\
&=&E_{s,x}|\varphi(X^{s,x}_T)|^2 +2E_{s,x}\int_t^T(f_u
u_\varepsilon+u_\varepsilon r^\varepsilon)(\theta,X_\theta)\,d\theta\\
&\le&E_{s,x}|\varphi(X^{s,x}_T)|^2
+E_{s,x}\int_t^T|f_u(\theta,X_\theta)|^2\,d\theta
+E_{s,x}\int_t^T|u_{\varepsilon}(\theta,X_\theta)|^2\,d\theta\\
&&+|I^1_{\varepsilon}|
+2E_{s,x}\int_t^Tr^{\varepsilon}h^+(\theta,X_\theta)\,d\theta
\end{eqnarray*}
for every $t\in[s,T]$. Hence, by Gronwall's lemma,
\begin{equation}
\label{eq3.20} I^2_{\varepsilon}\le
CE_{s,x}\left(|\varphi(X^{s,x}_T)|^2
+\int_s^T|f_u(t,X_t)|^2\,dt+|I^1_{\varepsilon}|
+\int_s^T(r^\varepsilon h^+)(t,X_t)\,dt\right).
\end{equation}
For any $\alpha>0$,
\[
E_{s,x}\int_s^T(r^\varepsilon h^+)(t,X_t)\,dt \le
E_{s,x}\left(\alpha \sup_{s\leq t\leq
T}|h^{+}(t,X_t)|^2+\alpha^{-1}|\int_s^T
r^\varepsilon(t,X_t)\,dt|^{2}\right)
\]
and, by (\ref{gwiazdka}),
\begin{equation}
\label{eq3.21} E_{s,x}|\int_s^T r^{\varepsilon}(t,X_t)\,dt|^{2}
\le E_{s,x}|\varphi(X^{s,x}_T)|^2
+E_{s,x}\int_s^T|f_u(t,X_t)|^2\,dt +I^2_{\varepsilon}.
\end{equation}
Hence, choosing a sufficiently large $\alpha$ we see from
(\ref{eq3.20}) that
\begin{equation}
\label{eq3.022} I^2_{\varepsilon}\le
CE_{s,x}\left(|\varphi(X^{s,x}_T)|^2 +\int_s^T|f_u(t,X_t)|^2\,dt
+|I^1_{\varepsilon}| +\sup_{s\leq t\leq T}|h^{+}(t,X_t)|^2\right).
\end{equation}
Therefore,  combining (\ref{eq3.21}) with (\ref{eq3.31}),
(\ref{eq3.022}) and  using Fatou's lemma we conclude that
$E_{s,x}(\int_s^Tr(t,X_t)\,dt)^2<\infty$. Finally, by
(\ref{gwiazdka}) and It\^o's formula, for any
$\varepsilon_1,\varepsilon_2>0$ and $\alpha>0$ we have
\begin{eqnarray*}
&&E_{s,x}|(u_{\varepsilon_1}-u_{\varepsilon_2})(t,X_t)|^2
+E_{s,x}\int^T_s
|\sigma\nabla(u_{\varepsilon_1}-u_{\varepsilon_2})(t,X_t)|^2\,dt\\
&&\quad\le CE_{s,x}\int^T_s(r^{\varepsilon_1}-r^{\varepsilon_2})
(u_{\varepsilon_1}-u_{\varepsilon_2})(t,X_t)\,dt\\
&&\quad\le CE_{s,x}\sup_{s\le t\le T}
|(u_{\varepsilon_1}-u_{\varepsilon_2})(t,X_t)|
\int^T_s|(r^{\varepsilon_1}-r^{\varepsilon_2})(t,X_t)|\,dt\\
&&\quad\le\alpha^{-1}CE_{s,x}\sup_{s\le t\le T}
|(u_{\varepsilon_1}-u_{\varepsilon_2})(t,X_t)|^2 +\alpha C
E_{s,x}\left(\int^T_s
|(r^{\varepsilon_1}-r^{\varepsilon_2})(t,X_t)|\,dt\right)^2.
\end{eqnarray*}
Hence, using the  Burkholder-Davis-Gundy inequality we obtain the
estimate
\begin{eqnarray*}
&&E_{s,x}\sup_{s\le t\le T}
|(u_{\varepsilon_1}-u_{\varepsilon_2})(t,X_t)|^2 +E_{s,x}\int^T_s
|\sigma\nabla(u_{\varepsilon_1}-u_{\varepsilon_2})(t,X_t)|^2\,dt\\
&&\qquad\le CE_{s,x}\left(\int^T_s
|(r^{\varepsilon_1}-r^{\varepsilon_2})(t,X_t)|\,dt\right)^2
\end{eqnarray*}
with $C$ not depending on $\varepsilon_1,\varepsilon_2$. Therefore
letting $\varepsilon\downarrow0$ in (\ref{gwiazdka}) we see that
the triple $(u(t,X_{t}),\sigma\nabla
u(t,X_{t}),\int_{0}^{t}r(t,X_{t}^{s,x})\,dt)$, $t\in[s,T]$, is a
solution of RBSDE$(\varphi,f,h)$ which in view of uniqueness
completes the proof.
\end{dow}

\begin{lm}\label{uw1}
If $\tilde{u}$ is a solution of PDE$(\varphi,f)$ and $(u,\mu)$ is a
solution of OP$(\varphi,f,h\vee\tilde{u})$ then $(u,\mu)$ is a
solution of OP$(\varphi,f,h)$.
\end{lm}
\begin{dow}
Let $(u,\mu)$ be a solution of OP$(\varphi,f,h\vee\tilde{u})$. Then
$u\geq h\vee\tilde{u}\geq h$. Moreover, by comparison results, for
any solution $(u_1,\mu_1)$ of  OP$(\varphi,f,h_1)$ with some $h_{1}$ we have
$u_1\geq\tilde{u}$. Hence $\mu_1\mathbf{1}_{\{h_1\leq
\tilde{u}\}}=0$, and consequently,
\[
\int_{Q_T}(u-h)\,d\mu=\int_{\{h\leq
\tilde{u}\}}(u-h)\,d\mu+\int_{\{h>\tilde{u}\}}(u-h)\,d\mu
=\int_{\{h>\tilde{u}\}}(u-(h\vee\tilde{u}))\,d\mu=0,
\]
which proves the lemma.
\end{dow}

\begin{lm}
\label{lm3.2}
Let $\varphi\in \mathbb{L}_{2,\varrho}(Q_{T})$, $g\in
\mathbb{L}_{p,q,\varrho}(Q_{T})$. Then
\[
E_{s,x}|\varphi(X_{T})|^{2}\le
C\varrho^{-2}(x)(T-s)^{-d/2}\|\varphi\|_{2,\varrho}^{2}
\]
and
\[
E_{s,x}\int_{s}^{T}|g(t,X_{t})|^{2}\,dt\le
C\varrho^{-2}(x)\|g\|^{2}_{p,q,\varrho}.
\]
\begin{dow}
Both inequalities follows form Aronson's estimates, because
\begin{eqnarray*}
\int_{\mathbb{R}^d}|\varphi(y)|^2p(s,x,T,y)\,dy
&\le&C\varrho^{-2}(x)\int_{\mathbb{R}^d}|\varphi(y)|^2|\varrho(y)|^2
\frac{p(s,x,T,y)}{|\varrho(y-x)|^2}\,dy \\
&\le& C\varrho^{-2}(x)(T-s)^{-d/2}\|\varphi\|^{2}_{2,\varrho}
\end{eqnarray*}
and, by H\"older's inequality,
\begin{eqnarray*}
&&\int_s^T\!\!\int_{\mathbb{R}^d}|g(t,y)|^2p(s,x,t,y)\,dt\,dy\\
&&\qquad\le C\varrho^{-2}(x)\int_s^T\!\!
\int_{\mathbb{R}^d}|g(t,y)|^2
\varrho^2(y) p(s,x,t,y)\varrho^{-2}(y-x)\,dt\,dy\\
&&\qquad \le C\varrho^{-2}(x)\|g\|^{2}_{p,q,\varrho}
\|p(0,0,\cdot,\cdot)\varrho^{-2}\|_{(p/2)^{*},(q/2)^{*}},
\end{eqnarray*}
which is finite by Aronson's estimate.
\end{dow}
\end{lm}

\begin{lm}
\label{lm3.3} If $\varphi\in\mathbb{L}_{2,\varrho}(\mathbb{R}^d)$,
$g \in \mathbb{L}_{2,\varrho}(Q_{T})$ for some $\varrho\in W$ and
(\ref{eq2.05}) is satisfied for every $(s,x)\in
[0,T)\times\mathbb{R}^{d}$, then for every
$K\subset\subset[0,T)\times\mathbb{R}^d$
\[
\sup_{(s,x)\in K}
E_{s,x}|(\varphi-\varphi_{n})(X_{T})|^{2}\rightarrow0
\]
and
\[
\sup_{(s,x)\in
K}E_{s,x}\int_{s}^{T}|(g-g_{n})(t,X_{t})|^{2}\,dt\rightarrow0
\]
as $n\rightarrow\infty$, where
$\varphi_n=\varphi\mathbf{1}_{B(0,n)}$, $g_n=g\mathbf{1}_{B(0,n)}$.
\end{lm}
\begin{dow}
The first assertion follows from Lemma \ref{lm3.2}. To prove the
second, let us choose $R>0$ such that $K\subset [0,T)\times
B(0,R)$ and $x\in B(0,R)$. Then for $n\ge2R$ we have
\begin{eqnarray*}
&&\int_s^T\!\!\int_{\mathbb{R}^d}|(g_n-g)(t,y)|^2p(s,x,t,y)\,dt\,dy\\
&&\qquad=\int_s^T\!\! \int_{B^{c}(0,2R)}|(g_n-g)(t,y)|^2
\varrho^2(y)p(s,x,t,y)\varrho^{-2}(y)\,dt\,dy\\
&&\qquad \le C\varrho^{-2}(x)\int_s^T\!\!\int_{B^{c}(0,2R)}
|(g_n-g)(t,y)|^2\varrho^{2}(y)\psi(s,x,t,y)\,dt\,dy,
\end{eqnarray*}
where $\psi(s,x,t,y)=(t-s)^{-d/2}\exp(-\frac{|y-x|^{2}}
{C(t-s)})(1+|y-x|^2)^{2\alpha}$. Since $\psi$ is bounded for $0\le
s<t\le T $, $|x-y|>R$  we see that
\begin{eqnarray}
\label{eq3.26} E_{s,x}\int_{s}^{T}|(g-g_{n})(t,X_{t})|^{2}\,dt
&=&\int_s^T\!\!\int_{\mathbb{R}^d}
|(g_n-g)(t,y)|^2p(s,x,t,y)\,dt\,dy\nonumber\\
&\le& C\varrho^{-2}(x)\|g_{n}-g\|^2_{2,\varrho,T}
\end{eqnarray}
for $(s,x)\in K$, $n\ge2R$, which completes the proof.
\end{dow}

\begin{tw}
\label{tw3.8} Let assumptions (H1)-(H3) hold with
$\varphi\in\mathbb{L}_{2,\varrho}(\mathbb{R}^{d})$,
$g\in\mathbb{L}_{2,\varrho}(Q_{T})$, where $\varrho\in W$ and
$\varrho<1$. Moreover, assume that  $g$ satisfies (\ref{n2.2}),
$h\in C(Q_{T})$ and $h\le c\bar\varrho^{-1}$ for some $c>0$ and
$\bar\varrho\in W$ such that
$\bar\varrho^{-1}\in\mathbb{L}_{2,\varrho}(\mathbb{R}^d)$. Then
there exists a unique solution $(u,\mu)$ of OP$(\varphi,f,h)$ such
that $u\in C([0,T)\times \mathbb{R}^{d})\cap
W_{2,\varrho}^{0,1}(Q_{T})$ and (\ref{eq3.02}), (\ref{eq3.04})
hold for each $(s,x)\in[0,T)\times \mathbb{R}^{d}$.
\end{tw}
\begin{dow}
We divide the proof into two steps: the case of linear and
semilinear equation.\\
{\em Step 1}. We first assume that $f=f(t,x),(t,x)\in Q_T$
satisfies (\ref{n2.2}) with $g$ replaced by $f$. Suppose that
$h(x)\le c\bar\varrho^{-1}$, where
$\bar\varrho(x)=(1+|x|^2)^{-\beta}$ for some $c,\beta>0$. Set
$\varphi_n=\mathbf{1}_{B(0,n)}\varphi$, $f_n=\mathbf{1}_{B(0,n)}f$
and consider a sequence $\{h_n\}\subset W^{1,1}_2(Q_T)$ such that
$h_n\le 2c\bar\varrho^{-1}$, $n\in\mathbb{N}$, and $h_n\rightarrow
h$ uniformly in compact subsets of $Q_T$. By Theorem \ref{tw3.4},
for each $n\in\mathbb{N}$ there is a unique solution $(u_n,\mu_n)$
of OP$(\varphi_n,f_n,h_n)$, and moreover,
\[
(u_n(t,X_t),\sigma\nabla
u_n(t,X_t))=(Y^{s,x,n}_{t},Z^{s,x,n}_t),\quad P_{s,x}\mbox{-}a.s.
\]
and
\[
\int_s^T\!\!\int_{\mathbb{R}^d}\xi(t,y)p(s,x,t,y)\,d\mu_n(t,y)
=E_{s,x}\int_s^T\xi(t,X_t)\,dK^{s,x,n}_t
\]
for all $\xi\in C_b(Q_T)$, where
$(Y^{s,x,n}_t,Z^{s,x,n}_t,K^{s,x,n}_t)$ is a solution of
RBSDE$(\varphi_n,f_n,h_n)$. By Lemma \ref{oszacowanie przyrostu},
\begin{eqnarray}
\label{nierownosc1} &&E_{s,x}\sup_{s\le t\le T}|(u_n-u_m)(t,X_t)|^2
+E_{s,x}\sup_{s\le
t\le T}|K^{s,x,n}_t-K^{s,x,m}_t|^2\nonumber\\
&&\qquad+E_{s,x}\int_s^T|\sigma\nabla(u_n-u_m)(t,X_t)|^2\,dt \nonumber\\
&&\quad\le E_{s,x}\sup_{s\le t\le T}|(h_n-h_m)(t,X_t)|^2
+E_{s,x}|(\varphi_n-\varphi)(X_T)|^2 \nonumber\\
&&\qquad+E_{s,x}\int_s^T |(f_n-f_m)(\theta,X_\theta)|^2\,d\theta.
\end{eqnarray}
From this and Theorem \ref{tw2.2} we deduce that
\begin{eqnarray}
\label{eq3.22} &&\|(u_{n}-u_m)(s)\|^{2}_{2,\varrho}
+\|\nabla(u_{n}-u_{m})\|^{2}_{2,\varrho^,T} \nonumber\\
&&\quad\le C\left(\int_{\mathbb{R}^d} E_{s,x}\sup_{s\le t\le T}|
(h_n-h_m)(t,X_t)|^2 \varrho^{2}(x) \,dx
+\|\varphi_{n}-\varphi_{m}\|^{2}_{2,\varrho}\right.\nonumber\\
&&\qquad\qquad\left.+\|f_{n}-f_{m}\|^{2}_{2,\varrho,T}\right).
\end{eqnarray}
Using Theorem \ref{tw2.2} we also get
\begin{equation}
\label{eq3.23} \sup_{0\le t\le T}
\|(u_{n}-u_m)(t)\|^{2}_{2,\varrho}\le C\int_{\mathbb{R}^d}\sup_{0\le
t\le T} E_{0,x}|(u_n-u_m)(t,X_t)|^2 \varrho^2(x)\,dx.
\end{equation}
Due to Lemma \ref{uw1}, without loss of generality we may assume
that $h_n\geq\tilde{u}_n$, where $\tilde{u}_n$ is a solution of
PDE$(\varphi_n,f_n)$. From comparison theorem (see \cite{EKPPQ}) we know that
$\underline{u}\le\tilde{u}_n$, where $\underline{u}$ is a
continuous solution of PDE$(-|\varphi|,-|f|)$, and that
$\underline{u}_n\searrow \underline{u}$, where $\underline{u}_n$
is a continuous solution of PDE$(-|\varphi_n|,-|f_n|)$. Since
\[
E_{s,x}\sup_{s\le t\le T}|(\underline{u}(t,X_t))^{+}|^2 \le
E_{s,x}\sup_{s\le t\le T}|(\underline{u}_n(t,X_t))^{+}|^2
\]
and
\begin{eqnarray*}
E_{s,x}\sup_{s\le t\le T}|(\underline{u}(t,X_t))^{-}|^2
&=&E_{s,x}\lim_{n\rightarrow\infty}
\sup_{s\le t\le T}|(\underline{u}_n(t,X_t))^{-}|^2\\
&\le&\liminf_{n\rightarrow\infty}E_{s,x}\sup_{s\le t\le T}|
(\underline{u}_n(t,X_t))^{-}|^2,
\end{eqnarray*}
from a priori estimates for solutions of
BSDE$(-|\varphi_{n}|,-|f_{n}|)$ (see \cite{A.Roz.BSDE}) we get
\begin{equation}
\label{gwiazdka1} E_{s,x}\sup_{s\le t\le
T}|\underline{u}(t,X_t)|^2 \le CE_{s,x}\left(|\varphi(X_T)|^2
+\int_s^T|g(t, X_t)|^2\,dt\right).
\end{equation}
Since
$|h_{n}(t,X_{t})|\le|\underline{u}(t,X_{t})|+2(1+|X_t|^2)^{\beta}$
and $\{h_n\}$ converges uniformly in compact subsets of $Q_T$, using
(\ref{gwiazdka1}), Proposition \ref{doob} and Lemma \ref{lm3.3} we
conclude that the right-hand side of (\ref{nierownosc1}) converges
to zero as $n,m\rightarrow\infty$.
From this it follows that there is $u$ such that $u_n\rightarrow
u$ pointwise in $[0,T)\times\mathbb{R}^d$. Moreover, using
(\ref{eq3.22}), (\ref{eq3.23}) and arguing as in the proof of
convergence of the right-hand side of (\ref{eq3.010}) we conclude
that $u_n\rightarrow u$ in $W^{0,1}_{2,\varrho}(Q_T)$ and
$u_n\rightarrow u$ in
$C([0,T],\mathbb{L}_{2,\varrho}(\mathbb{R}^d))$. By the definition
of solution of OP$(\varphi_n,f_n,h_n)$,
\begin{eqnarray}
\label{eq3.30} &&\int_t^T\langle u_n(s),\frac{\partial
\eta}{\partial s}(s)\rangle_2\,ds+\frac12\int_t^T\langle a(s)\nabla
u_n(s),\nabla\eta(s)\rangle_2\,ds\nonumber\\
&&\quad=\int_t^T \langle f_n(s),\eta(s)\rangle_2\,ds
+\int_t^T\!\!\int_{\mathbb{R}^d}\eta\, d\mu_n
+\langle\varphi_n,\eta(T)\rangle_2-\langle u_n(t),\eta(t)\rangle_2
\end{eqnarray}
for any $\eta\in W^{1,1}_{2,0}(Q_T)$. Therefore, if
$K\subset\subset[0,T)\times \mathbb{R}^d$, then choosing $\eta\in
W^{1,1}_2(Q_T)\cap C_0(Q_T)$ such that $\eta\equiv1$ on $K$ and
$0\le \eta\le 1$ we deduce from (\ref{eq3.30}) and Proposition
\ref{stw2.30} applied to $(u_n,\mu_n)$ and
$p\equiv2\bar\varrho^{-1}$ that $\sup_{n\ge1}\mu_n(K)<\infty$.
Thus, $\{\mu_n\}$ is tight. Taking a subsequence if necessary we
may assume that $\mu_n\Rightarrow\mu$, where $\mu$ is a Radon
measure on $Q_T$. Using  arguments similar to those in the proof
of {\em Step 3}\, of Theorem \ref{tw3.4} shows that
$\mu(\{t\}\times\mathbb{R}^d)=0$ for every $t\in[0,T]$. Now we
will show that $u$ is continuous on $[0,T)\times\mathbb{R}^d$. By
Lemma \ref{oszacowanie przyrostu}, for any $0<\delta<T/2$ and
$K\subset\subset[0,T-2\delta]\times\mathbb{R}^d$,
\begin{eqnarray*}
&&\sup_{(s,x)\in K}|(u_n-u_m)(s,x)|^2\le\sup_{(s,x)\in
K}E_{s,x}\sup_{s\le t\le
T-\delta}|(h_n-h_m)(t,X_t)|^2\\
&&\qquad+\sup_{(s,x)\in K}\int_s^T|(f_n-f_m)(t,X_t)|^2\,dt
+\sup_{(s,x)\in K}
E_{s,x}|Y^{s,x,n}_{T-\delta}-Y^{s,x,m}_{T-\delta}|^2.
\end{eqnarray*}
Since $h_n\rightarrow h$ uniformly in compact subsets of $Q_T$,
using Theorem \ref{tw2.1}, Proposition \ref{doob} and Lemma
\ref{lm3.4} we conclude that the first term on the right-hand of
above inequality converges to zero as $n,m\rightarrow\infty$.
Convergence of the second term follows from Lemma \ref{lm3.3} and
the third from the inequality
\[
\sup_{(s,x)\in K}
E_{s,x}|Y^{s,x,n}_{T-\delta}-Y^{s,x,m}_{T-\delta}|^2
\le\sup_{(s,x)\in K}C\varrho^{-2}(x)\delta^{-d/2}
\|(u_n-u_m)(T-\delta)\|_{2,\varrho}
\]
which is a consequence of Lemma \ref{lm3.2}. Thus, $u\in
C([0,T)\times\mathbb{R}^d)$. The pair $(u,\mu)$ is a weak solution
of OP$(\varphi,f,h)$, because letting $n\rightarrow\infty$ in
(\ref{eq3.30}) we get (\ref{eq2.1}), and similarly, by passing to
the limit we show that conditions (b), (d) of the definition of a
solution are satisfied. Finally, (\ref{eq3.02}), (\ref{eq3.04}) we
show as in the proof of Theorem \ref{tw3.4}.
\medskip\\
{\em Step 2.} We consider the general semilinear case. For
$\gamma>0$ to be determined later let $V({\gamma})$ denote the
Banach space consisting of elements $u$ of
$W^{0,1}_{2,\varrho}(Q_T)\cap
C([0,T],\mathbb{L}_{2,\varrho}(\mathbb{R}^d))$ equipped with the
norm
\[
\|u\|_{V(\gamma)}^{2}=\sup_{0\le s\le T}
\|u_\gamma(s)\|^2_{2,\varrho}+\|u_\gamma\|^2_{2,\varrho,T}
+\frac{\lambda}{2}\|\nabla u_\gamma\|^2_{2,\varrho,T},
\]
where $u_{\gamma}(s,x)=e^{s\gamma/2}u(s,x)$. Write
$K_{n}=[0,T-1/n]\times B(0,n)$. By $\mathcal{W}(\gamma)$ we denote
the Fr\'echet space of elements of $W^{0,1}_2(Q_T)$ such that
$\sup_{(s,x)\in K_{n}}\|u\|_{\gamma,s,x}< \infty$ for all $n \in
\mathbb{N}$ equipped with the $F$-norm
\[
\|u\|_{\mathcal{W}(\gamma)}=\sum_{n=0}^{\infty}\frac{1\wedge
\sup_{(s,x)\in K_{n}} \|u\|_{\gamma,s,x}}{2^{n}}\,,
\]
where $\|u\|^2_{\gamma,s,x}
=\int_s^T\!\!\int_{\mathbb{R}^d}e^{\gamma t}(|u(t,y)|^2
+\lambda|\nabla u(t,y)|^2)p(s,x,t,y)\,dt\,dy$, and by $B$ the
Fr\'echet  space of elements of $C([0,T)\times\mathbb{R}^{d})$ with
the $F$-norm
\[
\|u\|_{B} =\sum_{n=1}^{\infty}\frac{1\wedge
\|u\mathbf{1}_{K_n}\|_{\infty}}{2^{n}}\,.
\]
Finally, let $\mathcal{M}_{\gamma}$ denote the Fr\'echet space
$B\cap V(\gamma)\cap\mathcal{W}(\gamma)$ equipped with the
$F$-norm
$|\!|\!|u|\!|\!|_{\gamma}=\|u\|_{V(\gamma)}+\|u\|_{\mathcal{W}(\gamma)}
+\|u\|_{B}$. Now, define the mapping
$\Phi:\mathcal{M}_{\gamma}\rightarrow\mathcal{M}_{\gamma}$ by
putting $\Phi(v)$ to be the first component of the solution
$(u,\mu)$ of OP$(\varphi,f_{v},h)$. By {\em Step 1} the definition
of $\Phi$ is correct. We are going to show that $\Phi$ is
contractive on $\mathcal{M}_{\gamma}$. Let
$v_1,v_2\in\mathcal{M}_{\gamma}$ and let $(u_i,\mu_i)$, $i=1,2$,
be solutions of OP$(\varphi,f_{v_i},h)$. Set
$u=u_1-u_2=\Phi(v_1)-\Phi(v_2)$, $\mu=\mu_1-\mu_2$. By the
definition of a solution of the obstacle problem,
\begin{eqnarray*}
&&\langle u(t),\eta(t)\rangle_2+\int_t^T\langle u(s),\frac{\partial
\eta}{\partial s}(s)\rangle_2\,ds+\frac12\int_t^T\langle a(s)\nabla
u(s),\nabla\eta(s)\rangle_2\,ds\\
&&\qquad=\int_t^T\!\!\int_{\mathbb{R}^d}\eta d\mu+\int_t^T\langle
f_{v_1}(s)-f_{v_2}(s),\eta(s)\rangle_2\,ds.
\end{eqnarray*}
Putting $\eta(s)=e^{\gamma s}u(s)\varrho^2\xi_n^2$, where $\xi_n$
is defined as in the proof of Proposition \ref{stw2.30}, we obtain
\begin{eqnarray*}
&&e^{\gamma t}\langle u(t),u(t)\varrho^2\xi_n^2
\rangle_2+\gamma\int_t^Te^{\gamma s}\langle
u(s),u(s)\varrho^2\xi_n^2 \rangle_2\,ds\\
&&\qquad+\int_t^Te^{\gamma s}\langle u(s),\frac{\partial u}{\partial
s}(s)\varrho^2\xi_n^2\rangle\,ds +\frac12\int_t^Te^{\gamma s}
\langle a(s)\nabla
u(s),\nabla(u\varrho^2\xi_n^2)(s)\rangle_2\,ds\\&&
\quad=\int_t^T\!\!\int_{\mathbb{R}^d}u(s)\varrho^2\xi_n^2\,d\mu
+\int_t^Te^{\gamma s}\langle
f_{v_1}(s)-f_{v_2}(s),u(s)\varrho^2\xi^2_n\rangle\,ds.
\end{eqnarray*}
By the above and (\ref{eq3.03}),
\begin{eqnarray*}
&& e^{\gamma t}\|u(t)\varrho\xi_n\|^2_2 +\gamma\int_t^Te^{\gamma s}
\|u(s)\varrho\xi_n\|^2_2\,ds +\frac12\int_t^Te^{\gamma s}
\frac{d}{ds}\|u(s)\varrho\xi_n\|^2\,ds\\
&&+\frac12\int_t^Te^{\gamma s}\langle a(s)\nabla
u(s),\nabla(u\varrho^2\xi_n^2)(s)\rangle_2\,ds\le \int_t^Te^{\gamma
s} \langle f_{v_1}(s)-f_{v_2}(s),u(s)\varrho^2\xi^2_n\rangle\,ds.
\end{eqnarray*}
Consequently,
\begin{eqnarray*}
&& e^{\gamma t}\|u(t)\varrho\xi_n\|^2_2+\gamma\int_t^Te^{\gamma
s}\|u(s)\varrho\xi_n\|^2_2\,ds+\int_t^Te^{\gamma s}\langle
a(s)\nabla u(s),\nabla(u\varrho^2\xi_n^2)(s)\rangle_2\,ds\\
&&\qquad\le 2\int_t^Te^{\gamma s}\langle
f_{v_1}(s)-f_{v_2}(s),u(s)\varrho^2\xi^2_n\rangle_2\,ds.
\end{eqnarray*}
Letting $n\rightarrow\infty$ and performing computations similar to
that in the proof of Theorem \ref{jednoznacznosc} we get
\begin{eqnarray*}
&&e^{\gamma t}\|u(t)\varrho\|^2_2 +\gamma\int_t^Te^{\gamma s}
\|u(s)\varrho \|^2_2\,ds +\frac{\lambda}{2}\int_t^Te^{\gamma s}
\|\nabla u(s)\varrho\|^2_2\,ds\\
&&\quad\le 2\int_t^Te^{\gamma s} \langle
f_{v_1}(s)-f_{v_2}(s),u(s)\varrho^2\rangle_2\,ds
+\frac{\Lambda}{2\lambda}\int_t^Te^{\gamma s}
\|u(s)\varrho\|^2_2\,ds.
\end{eqnarray*}
The right-hand side of the above inequality may be estimated by
\begin{eqnarray*}
&&2L\int_t^Te^{\gamma s}
\|(v_1-v_2)(s)\varrho\|_2\|u(s)\varrho \|_2\,ds\\
&&\qquad+2L\Lambda\int_t^T e^{\gamma s}
\|\nabla(v_1-v_2)(s)\varrho\|_2 \|u(s)\varrho\|_2\,ds
+\frac{\Lambda}{2\lambda}\int_t^Te^{\gamma s}\|u(s)\varrho\|^2_2\,ds\\
&&\quad\le(4L^{2}+\frac{8\Lambda^{2}L^2}{\lambda}
+\frac{\Lambda}{2\lambda})\int_t^Te^{\gamma s}
\|u(s)\varrho\|^2_2ds+\frac14\int_t^Te^{\gamma
s}\|(v_1-v_2)(s)\varrho\|^2_2\,ds \\
&&\qquad+ \frac{1}{4}\int_t^T \frac{\lambda}{2}e^{\gamma s}
\|\nabla(v_1-v_2)(s)\varrho\|^2_2\,ds.
\end{eqnarray*}
Putting
$\gamma=1+4L^{2}+8\lambda^{-1}\Lambda^2L^2+(2\lambda)^{-1}\Lambda$
we see that
\begin{eqnarray}\label{gwiazdka2}
\|\Phi(v_1)-\Phi(v_2)\|_{V(\gamma)}\le2^{-1}\|v_1-v_2\|_{V(\gamma)}.
\end{eqnarray}
Let $(Y^{s,x,i},Z^{s,x,i},K^{s,x,i})$, $i=1,2$, denote a solution of
RBSDE$(\varphi,f_{v_i},h)$ and let $v=v_1-v_{2}$. We already know
that $(Y^{s,x,i}_t,Z^{s,x,i}_t)=(u_i(t,X_t),\sigma\nabla
u_i(t,X_t))$, $t\in[s,T]$. Therefore, since
$E_{s,x}\int_t^Te^{\gamma \theta}
v(\theta,X_\theta)\,d(K^{s,x,1}_\theta-K^{s,x,2}_\theta)\le0$ for
every $t\in[s,T]$, using It\^o's formula we have
\begin{eqnarray*}
&&E_{s,x}e^{\gamma t}|u(t,X_t)|^2
+E_{s,x}\int_t^Te^{\gamma\theta}(\gamma|u(\theta,X_\theta)|^2
+|\sigma\nabla u(\theta,X_\theta)|^2)\,d\theta\\
&&\qquad\le2E_{s,x}\int_t^Te^{\gamma \theta}
u(f_{v_1}-f_{v_{2}})(\theta,X_\theta)\,d\theta\\
&&\qquad\le 2LE_{s,x}\int_t^Te^{\gamma \theta}
u(\theta,X_\theta)(|v|+|\sigma\nabla v|)(\theta,X_\theta)d\theta\\
&&\qquad\le 8\lambda^{-1}\Lambda L^2\varepsilon
E_{s,x}\int_t^Te^{\gamma\theta}
|u(\theta,X_\theta)|^2\,d\theta\\
&&\qquad\quad +\varepsilon^{-1}E_{s,x}\int_t^Te^{\gamma \theta}
(|v|^2+\lambda|\nabla v|^2)(\theta,X_\theta)|\,d\theta.
\end{eqnarray*}
Putting $\gamma=1+8\Lambda\lambda^{-1}L^2\varepsilon$  with suitably
chosen $\varepsilon>0$ in a standard manner we deduce from the above
and the Burkholder-Davis-Gundy inequality that
\begin{eqnarray}
\label{dodatkowa-dwie gwiazdki} &&E_{s,x}\sup_{s\le t\le T}
e^{\gamma t}|u(t,X_t)|^2 +E_{s,x}\int_s^T\lambda e^{\gamma t}
(|u(t,X_{t})|^2+|\nabla u(t,X_{t})|^2)\,dt\nonumber\\
&&\qquad\le4^{-1}E_{s,x}\int_s^Te^{\gamma t} (|v|^2+\lambda|\nabla
v|^2)(t,X_t)\,dt.
\end{eqnarray}
From this we obtain
\[
\|\Phi(v_{1})-\Phi(v_{2})\|_{B}
+\|\Phi(v_{1})-\Phi(v_{2})\|_{\mathcal{W}(\gamma)}\leq
2^{-1}(\|v_{1}-v_{2}\|_{B} +\|v_{1}-v_{2}\|_{\mathcal{W}(\gamma)}),
\]
which when combined with (\ref{gwiazdka2}) shows that $\Phi$ is
contractive on $\mathcal{M}_{\gamma}$.  By Banach's principle,
$\Phi$ has a unique fixed point $u$. Clearly, the solution $(u,\mu)$
of OP$(\varphi,f_u,h)$ has the asserted properties.
\end{dow}
\medskip

One can prove Theorem \ref{tw3.8} by the method of stochastic
penalization used in the proof of Theorem \ref{tw3.4}. To apply
that method one should first generalize results of
\cite{A.Roz.BSDE} on representation of solutions of the Cauchy
problem proved for $\varphi\in\mathbb{L}_{2}(\mathbb{R}^d)$,
$g\in\mathbb{L}_{2}(Q_T)$ to the case
$\varphi\in\mathbb{L}_{2,\varrho}(\mathbb{R}^d)$,
$g\in\mathbb{L}_{2,\varrho}(Q_T)$ for some $\varrho\in W$. Since
detailed proof of such a generalization does not bring new ideas
and at the same time requires some efforts, we decided to present
a different approach. Note, however, that the adopted approach
uses some ideas from \cite{A.Roz.BSDE}.

\begin{wn}\label{wn3.1}
Let assumptions of Theorem \ref{tw3.8} hold. Define
$(u_{n},\mu_{n})$ as in Theorem \ref{tw3.4}. Then
\begin{enumerate}
\item[\rm(i)]
$u_{n}\uparrow u$ uniformly in compact subsets of $[0,T)\times
\mathbb{R}^{d}$, $u_{n}\rightarrow u$ in
$W^{0,1}_{2,\varrho}(Q_{T})\cap
C([0,T],\mathbb{L}_{2,\varrho}(Q_{T}))$,
\item[\rm(ii)] $\mu_{n}\Rightarrow\mu$,
$\mu_{n}\rightarrow\mu$ in $(W_{2,\varrho}^{1,1}(Q_{T}))^{*}$,
$\mu_n\Rightsx\mu$ for every $(s,x)\in[0,T)\times\mathbb{R}^d$.
\end{enumerate}
\end{wn}
\begin{dow}
Follows from Theorems \ref{tw2.2}, \ref{tw2.4} and \ref{tw3.8}.
\end{dow}
\medskip

Let us remark that Corollaries \ref{wn3.5}, \ref{wn3.6} hold also
under the assumptions of Theorem \ref{tw3.8}. The proof of
Corollary \ref{wn3.5} runs as before. In the proof of Corollary
\ref{wn3.6} the main difference consists in the fact that instead
of boundedness of $\{u_{\varepsilon}\}$ in $W^{0,1}_2(Q_T)$ (see
(\ref{eq3.28})) we have to prove its boundedness in
$W^{0,1}_{2,\varrho}(Q_T)$. The last assertion one can show using
arguments from the proof of Proposition \ref{stw2.30}.

\begin{wn}
Let assumptions of Theorem \ref{tw3.4} or Theorem \ref{tw3.8} hold
and let $(u,\mu)$ be a solution of OP$(\varphi,f,h)$.
\begin{enumerate}
\item[\rm(i)] If  $g\in\mathbb{L}_{p,q,\varrho}(Q_{T})$ then
\[
|u(s,x)|+\|u\|_{\mathcal{W}_{2}(s,x,T)}\le C
{\varrho}^{-1}(x)(1+(T-s)^{-d/2} \|\varphi\|_{2,\varrho}^{2}
+\|f\|_{p,q,\varrho}^{2})^{1/2}.
\]
\item[\rm(ii)] If $|\varphi|\le c\varrho^{-1}$ for some
$c>0$, $\varrho\in W$  (i.e. $\varphi$ satisfies the polynomial
growth condition) then
\[
|u(s,x)|+\|u\|_{\mathcal{W}_{2}(s,x,T)}\le C
{\varrho^{-1}(x)}(1+\|f\|_{p,q,\varrho}^{2}))^{1/2}.
\]
\item[\rm(iii)]If $|\varphi|+|g|\le
c{\varrho}^{-1}$ for some $c>0$, $\varrho\in W$ then
\[
|u(s,x)|+\|u\|_{\mathcal{W}_{2}(s,x,T)}\le C {\varrho}^{-1}(x).
\]
\end{enumerate}
\end{wn}
\begin{dow}
It follows from (\ref{eq2.4}), Theorem \ref{tw3.8}, Proposition
\ref{doob} and  Lemma \ref{lm3.2}.
\end{dow}
\medskip

It is known that the obstacle problem (\ref{eq1.3}) with
non-divergent operator $A_t$ appear as the Hamilton-Jacobi-Bellman
equation for an optimal stopping time problem (see
\cite{BensoussanLions}) and that value functions of that stopping
problem is given by the first component of a solution of an RBSDEs
with forward driving processes associated with $A_t$ (see
\cite{EKPPQ}). It is worth noting that similar relations hold for
divergence form operators.

\begin{wn}
Let assumptions of Theorem \ref{tw3.4} or Theorem \ref{tw3.8} hold
and let $(u,\mu)$ be a solution of OP$(\varphi,f,h)$. Then for
each $t\in[s,T]$,
\begin{eqnarray*}
&&u(s,x)=\sup_{\tau\in \mathcal{T}^s_{t}} E_{s,x}(\int_{t}^{\tau}
f(\theta,X_{\theta},u(\theta,X_{\theta}),\sigma\nabla
u(\theta,X_{\theta}))\,d\theta\\
&&\qquad\qquad\qquad\qquad\qquad +h(\tau,X_{\tau})
\mathbf{1}_{\tau<T}+\varphi(X_{T})\mathbf{1}_{\tau=T}|\mathcal{G}^s_t),
\end{eqnarray*}
where $\mathcal{T}^s_{t}=\{\tau\in \mathcal{T}^s: t\le \tau\le
T\}$ and $\mathcal{T}^s$ denote the set of all
$\{\mathcal{G}^s_t\}$-stopping times.
\end{wn}
\begin{dow}
Let $\tau\in \mathcal{T}_{t}^{s}$. By (\ref{eq3.02}) and the
definition and a priori estimate for a solution of
RBSDE$(\varphi,f,h)$ we have
\begin{eqnarray*}
&&u(t,X_{t})=E_{s,x}(\int_{t}^{\tau}
f(\theta,X_{\theta},u(\theta,X_{\theta}),\sigma\nabla
u(\theta,X_{\theta}))\,d\theta+u(\tau,X_{\tau})
+K^{s,x}_{\tau}-K^{s,x}_{t}|\mathcal{G}^s_t)\\
&&\quad\geq
E_{s,x}(\int_{t}^{\tau}f(\theta,X_{\theta},u(\theta,X_{\theta}),\sigma\nabla
u(\theta,X_{\theta}))\,d\theta+h(\tau,X_{\tau})
\mathbf{1}_{\tau<T}+\varphi(X_{T})\mathbf{1}_{\tau=T}|\mathcal{G}^s_t).
\end{eqnarray*}
Let us define the optimal control by putting
$D_{t}=\inf\{t\le\theta\le T:
u(\theta,X_{\theta})=h(\theta,X_{\theta})\}$. Since  $K^{s,x}$ is
continuous and
$\int_{0}^{T}(u(t,X_{t})-h(t,X_{t}))\,dK_{t}^{s,x}=0$
$P_{s,x}$-a.s., it follows that $K^{s,x}_{D_{t}}-K^{s,x}_{t}=0$
$P_{s,x}$-a.s., which proves the corollary.
\end{dow}
\medskip

The next theorem provides a probabilistic formula for the minimal
weak solution of the variational inequality associated with
(\ref{eq1.3}).

\begin{tw}
\label{th3.16} Assume that (H1)--(H3) hold with $\varphi\in
\mathbb{L}_{2,\varrho}(\mathbb{R}^{d})$, $g\in
\mathbb{L}_{2,\varrho}(Q_{T})$ for some $\varrho\in W$ and with
$h\in C(Q_{T})$ satisfying the polynomial growth condition. Then
there exists a version $u$ of minimal weak solution of
OP$(\varphi,f,h)$ in the variational sense such that if
(\ref{eq2.05}) is satisfied for some $(s,x)\in [0,T)\times
\mathbb{R}^{d}$ then
\begin{equation}
\label{eq3.35} (Y_{t}^{s,x},Z_{t}^{s,x})=(u(t,X_{t}),\sigma\nabla
u(t,X_{t})),\quad t\in[s,T],\quad P_{s,x}\mbox{-}a.s.
\end{equation}
\end{tw}
\begin{dow}
By \cite[Theorem 4.1.6]{BensoussanLions} there exists the minimal
weak solution $\bar u$ of OP$(\varphi,f,h)$ in the variational
sense. Repeating arguments from the proof of Proposition
\ref{stw2.3} we show that there is a version $u$ of a weak
solution of the linear OP$(\varphi,f_{\bar u},h)$ in the
variational sense such that (\ref{eq3.35}) holds if (\ref{eq2.05})
is satisfied. Since $\|g\|_{2,\varrho}< \infty$, it follows that
(\ref{eq2.05}) is satisfied for a.e. $(s,x)\in [0,T)\times
\mathbb{R}^{d}$ (see remark following the proof of Theorem
\ref{tw2.4}). Therefore, by Theorems \ref{tw2.2} and  \ref{tw2.4},
$u$ is a limit in $W^{0,1}_{2,\varrho}(Q_T)$ of the penalizing
sequence defined by (\ref{eq3.2}), and hence (see the proof of
\cite[Theorem 4.1.6]{BensoussanLions}), $u$ is a minimal weak
solution of OP$(\varphi,f_{\bar u},h)$ in the variational sense.
Since the minimal solution is unique, $u=\bar u$, and the proof is
complete.
\end{dow}

\end{document}